\documentclass[11pt,reqno,a4paper]{amsart}
\usepackage{hyperref}
\usepackage{enumerate}
\usepackage[left=0.7in, right=0.7in, top= 0.8in, bottom=1.62in]{geometry}
\usepackage{amsmath,amssymb,amsthm,amsfonts}
\usepackage{graphicx}
\usepackage{tikz}

\usetikzlibrary{decorations.pathreplacing}
\usetikzlibrary{decorations.markings}
\usetikzlibrary{calc} 
\usepackage{subcaption} 

\usetikzlibrary{graphs, graphs.standard}

\usetikzlibrary{shapes}
\usetikzlibrary{plotmarks}
\usetikzlibrary{arrows}
\usetikzlibrary{positioning}
\theoremstyle{definition}
\newtheorem{defn}{Definition}[section]
\newtheorem{fact}[defn]{Fact}
\newtheorem{prop}[defn]{Proposition}
\newtheorem{thm}[defn]{Theorem}

\newtheorem{lem}[defn]{Lemma}

\newtheorem{remark}[defn]{Remark}

\newtheorem{claim}[defn]{Claim}

\title[Distinguishing invariants of Cayley Graphs and Group Actions]{Distinguishing invariants of Cayley Graphs and Group Actions$^{\ast}$}
\thanks{$^{\ast}$ An extended abstract of this work was presented at the Fourth Meru Combinatorics Conference 2026.}

\author{Amitayu Banerjee}
\address{E\"otv\"os Lor\'and University, Budapest, Hungary}
\email{banerjee.amitayu@gmail.com}

\date{}

\makeatletter
\@namedef{subjclassname@2020}{\textup{2020} Mathematics Subject Classification}
\makeatother
\subjclass[2020]{05C25, 05C15, 05E18, 20B25.}
\keywords{Cayley graph, group action, wreath product, Hamming graph, Cartesian powers of graphs, distinguishing number, distinguishing index, total distinguishing number}
\begin{document}

\begin{abstract} 
Let $\Gamma$ be a group acting on a graph $G=(V,E)$. The {\em total
distinguishing number} $D''_\Gamma(V\cup E)$ ({\em distinguishing index}
$D'_\Gamma(E)$) of the action is the least number of labels required
for a total labeling (edge labeling) of $G$ such that the only
elements of $\Gamma$ preserving the labeling belong to the
pointwise stabilizer of $V\cup E$ (of $E$). Building on the results of Aliyu (Ann. Commun. Math. 9, 2026), we obtain sharp bounds for these invariants
under free actions, wreath product
actions, and coordinatewise actions on Cartesian powers of graphs.
Alikhani and Soltani (J. Appl. Ind. Math. 15, 2021) proved that if $G=\operatorname{Cay}(H,S)$
is a normal Cayley graph, then its distinguishing number is at most
$1+D_{\operatorname{Aut}(H,S)}(S)$, where
$\operatorname{Aut}(H,S)=\{\alpha\in\operatorname{Aut}(H):
\alpha(S)=S\}$, for $S\subseteq H$.
We obtain analogous upper bounds for the distinguishing index and
total distinguishing number of normal Cayley graphs under certain
algebraic conditions.
Moreover, we present infinitely many graphs $G$ and subgroups $\Gamma \leq \operatorname{Aut}(G)$ with $D'_{\Gamma}(E) < D'_{\operatorname{Aut}(G)}(E)$ and $D''_{\Gamma}(V\cup E) < D''_{\operatorname{Aut}(G)}(V\cup E)$, and determine the distinguishing index for actions of $A_n$ on $K_n$ and $K_{1,n}$, of $A_3\times S_n$ on $K_{3,n}$, and of $S_n$ on Hamming graphs.
Finally, we show that for $n\geq4$, any faithful exactly $n$-transitive action on a set of size $m\geq n$ has distinguishing number at most $n+1$.
\end{abstract}

\maketitle
\section{Introduction} 
All graphs in this paper are finite, simple, and undirected, all groups are nontrivial and finite, and all sets are nonempty and finite. 
The notions of distinguishing number, distinguishing index, and total distinguishing number of graphs were introduced by Albertson and Collins \cite{AC1996}, Kalinowski and Pil\'{s}niak \cite{KP2015}, and Kalinowski et al. \cite{KPW2016}.

\begin{defn}\label{Definition 1.1}
Let $\mathrm{Aut}(G)$ be the automorphism group of a graph $G=(V,E)$, and $\mathrm{id}_{G}$ the identity automorphism of $G$.
Let $f$ be a vertex labeling (an edge labeling or a total labeling) of $G$. 
An automorphism $\varphi \in \mathrm{Aut}(G)$ is said to \emph{preserve} $f$ if:
\begin{enumerate}[(i)]
\item $f(x)=f(\varphi(x))$ for all $x\in V$, when $f$ is a vertex labeling;
\item $f(e)=f(\varphi(e))$ for all $e\in E$, when $f$ is an edge labeling;
\item $f(x)=f(\varphi(x))$ for all $x\in V\cup E$, when $f$ is a total labeling.
\end{enumerate}

A labeling $f$ is \emph{distinguishing} if it is preserved only by $\mathrm{id}_G$.
The distinguishing number $D(G)$ (distinguishing index $D'(G)$, total distinguishing number $D''(G)$) 
of $G$ is the least integer $d$ such that $G$ admits a distinguishing vertex (edge, total respectively) labeling with $d$ labels.
\end{defn}
\vspace{-2mm}

Tymoczko \cite{Tym2004} generalized the notion of distinguishing numbers to arbitrary group actions on sets. She defined the {\em distinguishing number of a group action} of a group $\Gamma$ on a set $X$, denoted by $D_{\Gamma}(X)$, as the minimum integer $k$ for which there exists a labeling $\varphi : X \to \{1,...,k\}$ preserved only by the elements of $\Gamma$ that belong to the pointwise stabilizer of $X$ in $\Gamma$.
Following this, the distinguishing number of group actions has been widely studied (cf. \cite{AS2021, AWTG2026, Ali2026, Chan2006, DP2023, WZ2009}, among others).
Alikhani and Soltani \cite[Theorem 2.2]{AS2021} proved that if $G = \mathrm{Cay}(H, S)$ is a normal Cayley graph, then $D(G)$ is at most $1+D_{\mathrm{Aut(H,S)}}(S)$
where $\mathrm{Aut}(H,S) = \{\alpha \in \mathrm{Aut}(H) \mid \alpha(S) = S\}$ for a group $H$ and a subset $S\subseteq H$.
Chan \cite{Chan2006}, Dadhwal and Pankaj \cite{DP2023}, Aliyu et al. \cite{AWTG2026}, and Aliyu \cite{Ali2026} studied the distinguishing number of wreath product and direct product actions. 
Recently, Aliyu et al. \cite{AWTG2026} and Aliyu \cite{Ali2026} proved the following results:

\begin{itemize}
    \item If a group $\Gamma$ acts on a set $V$ with a base of size $n$, then $D_{\Gamma}(V)\le n+1$.
    \item If $\Gamma$ acts freely on a set $V$, then $D_{\Gamma}(V) \leq 2$.
    \item If $\Gamma$ acts transitively on a set $V$, then for every integer $k \geq 1$, $D_{\Gamma}(V^k) \leq D_{\Gamma}(V) +k-1$.
    
    \item If $\Gamma$ acts faithfully on a set $V$, and $W = \Gamma \wr S_m$ act in the natural product action on $V^m$, then
$D_{W}(V^m) \leq D_{\Gamma}(V)+m-1.$
\end{itemize}
 
Inspired by the ideas of Tymoczko \cite{Tym2004}, we generalize the notion of distinguishing index and total distinguishing number from graph automorphisms to arbitrary group actions. 
Let $\Gamma$ be a group acting on a graph $G=(V,E)$. The {\em total distinguishing number of the action}, denoted $D''_{\Gamma}(V \cup E)$ ({\em distinguishing index of the action}, denoted $D'_{\Gamma}(E)$), is the minimum number of labels required for a total labeling (edge labeling) of $G$ so that the only elements of $\Gamma$ that preserve the labeling belong to the pointwise stabilizer of $V \cup E$ (pointwise stabilizer of $E$). 
If $\Gamma = \text{Aut}(G)$, these invariants coincide with the total distinguishing number and distinguishing index of $G$.
\vspace{-2mm}

\subsection{Results}
In Theorem \ref{Theorem 3.2}, 
we observe that if $G=\mathrm{Cay}(H,S)$ is a connected Cayley graph 
where $|H|,|S| \ge 2$, then $2 \le D''(G) \le 3$ and if $G$ is a graphical regular representation of $H$, then $D''(G)=2$.
Moreover, if $G$ is a normal Cayley graph, then similar to \cite[Theorem 2.2]{AS2021}, $D''(G)\leq 1+D''_{\mathrm{Aut(H,S)}}(S\cup E_{e}\cup\{e\})$ and $D'(G)\le 1 + D'_{\mathrm{Aut(H,S)}}(E_e)$ where $E_{e}=\{\{e,s\}:s\in S\}$ and $e$ is the identity of $H$, under certain algebraic conditions.
Let $S_{n}$ and $A_{n}$ denote the symmetric and alternating groups of degree $n$, and let $\mathbb{Z}_{n}$ denote the cyclic group of order $n$.
Let $C_n$, $P_n$, and $K_n$ denote the cycle, path, and complete graph, respectively, on $n$ vertices, and $K_{n,m}$ be a complete bipartite graph with bipartitions of size $m$ and $n$.
Let $k\geq 2$ be an integer.
Building on the research in \cite{AWTG2026} and \cite{Ali2026}, we prove the following results for a group $\Gamma$ acting on a graph $G=(V,E)$:

\begin{enumerate}
    \item (Theorem \ref{Theorem 4.4}) If the action of $\Gamma$ on $V$ has a base of size $n\geq 1$, then $D''_{\Gamma}(V\cup E)\leq n+1$, and if the induced action of $\Gamma$ on $E$ has a base of size $m\geq 1$, then $D'_{\Gamma}(E)\leq m+1$.
    The second bound is sharp for the actions of $S_p$ on $K_{1,p}$ for $p\geq 3$, and $A_q$ on $K_{1,q}$ for $q\geq 5$.
\vspace{2mm}
    \item (Propositions \ref{Proposition 4.5} and \ref{Proposition 4.6}) There exist infinitely many graphs $G=(V,E)$ and subgroups
    $\Gamma\leq\operatorname{Aut}(G)$ with
    $D'_{\Gamma}(E)<D'_{\operatorname{Aut}(G)}(E)$ and
    $D''_{\Gamma}(V\cup E)<D''_{\operatorname{Aut}(G)}(V\cup E)$.
\vspace{2mm}
    \item (Proposition \ref{Proposition 4.11}) If the action of $\Gamma$ on $G$ is free, then $D''_{\Gamma}(V\cup E)\leq 2$, and if no nontrivial element of $\Gamma$ of order $2$ swaps the endpoints of any edge, then $D'_{\Gamma}(E)\leq 2$.
    The bounds are sharp for the natural action of $\mathbb{Z}_p$ on $C_p$ for odd $p\geq 3$.
\vspace{2mm}
    \item (Theorem \ref{Theorem 5.4}) If the coordinatewise action of $\Gamma$ on $V^k$ induces an action of $\Gamma$ on the Cartesian power $G^k=(V^k,E_k)$, then
    $D''_{\Gamma}(V^k\cup E_k)\leq D''_{\Gamma}(V\cup E)$.
    The bound is sharp for the action of $\mathbb{Z}_2$ on $P_n$ for $n\geq 3$.
\end{enumerate}

We also obtain the following results for particular group actions:

\begin{enumerate}
    \item (Proposition \ref{Proposition 4.5}) For every $p\geq 2$,
    $D'_{A_p}(E(K_p))=2$.

    \item (Proposition \ref{Proposition 4.6}) For every $t\geq 5$,
    $
    D''_{A_t}(V(K_{1,t})\cup E(K_{1,t}))
    =\lceil(t-1)^{1/2}\rceil$.

    \item (Proposition \ref{Proposition 5.1}) If $n\geq 3$ and $S_n$ acts naturally on $K_n$ and coordinatewise on the Hamming graph $(K_n)^r=(V^r,E_r)$, then
    $
    D''_{S_n}(V^r\cup E_r)=D'_{S_n}(E_r)=2$
    for all $r\geq 2$.
\end{enumerate}

\vspace{-2mm}
\subsection{Main Results}
\begin{enumerate}
    \item (Theorem \ref{Theorem 4.7}) For every $n\geq 3$,
    $D'_{A_3\times S_n}\bigl(E(K_{3,n})\bigr)
=\left\lceil(n+1)^{1/3}\right\rceil$.

    \item (Theorem \ref{Theorem 4.14}) If $\Gamma$ is a group acting faithfully and exactly $n$-transitively on a set $X$ of size $m\geq n$, where $n\geq4$, then
$D_{\Gamma}(X)\leq n+1$.
    \item (Theorem \ref{Theorem 5.6}) If the action of a group $\Gamma$ on a graph $G=(V,E)$ is faithful, where $|V|\geq 2, |E|\geq 1$ and the natural action of the wreath product $W=\Gamma\wr S_k$ on $V^k$ induces an action of $W$ on $G^{k}=(V^{k},E_k)$,
then
$
D''_{W}(V^k\cup E_k)\le D''_{\Gamma}(V\cup E)+k-1
$ for any integer $k\geq 2$.
\end{enumerate}

\section{Preliminaries}

\begin{defn}\label{Definition 2.1}
Let $e$ be the identity of a group $\Gamma$. A group action of $\Gamma$ on a set $X$ is called:
\begin{enumerate}
    \item {\em free} if for any $x \in X$, the stabilizer of $x$ is trivial, i.e.,
    $
    \mathrm{Stab}_{\Gamma}(x)=\{g\in \Gamma: gx=x\}=\{e\}.
    $
    
    \item {\em faithful} if the stabilizer of $X$ is trivial, i.e.,
    $
    \mathrm{Stab}_{\Gamma}(X)=\{g\in \Gamma: gx=x (\forall x\in X)\}=\{e\}.
    $
    
    \item {\em $n$-transitive} if for any two $n$-tuples
    $(x_1,...,x_n)$ and $(y_1,...,y_n)$ of pairwise distinct
    elements of $X$, there exists at least one $g\in\Gamma$ such that
    $gx_i=y_i$ for every $i\in\{1,...,n\}$.

    \item {\em sharply $n$-transitive} if for any two 
    $n$-tuples $(x_1,...,x_n)$ and $(y_1,...,y_n)$ of pairwise
    distinct elements of $X$, there exists a unique $g\in\Gamma$
    such that $gx_i=y_i$ for every $i\in\{1,...,n\}$.
    \item {\em of base size $n$} if the smallest cardinality of a
subset $B\subseteq X$ such that $
\mathrm{Stab}_{\Gamma}(B)
=\{g\in\Gamma:gx=x\ (\forall x\in B)\}
=\{e\}
$ is $n$.
Such a subset $B$ is called a {\em base} for the action.
\end{enumerate}
\end{defn}
\vspace{-2mm}

If the action is free, then it is faithful, but the converse need not hold.
A sharply $n$-transitive action has base size exactly $n$, provided
$|X|\geq n+1$. A $1$-transitive action is called {\em transitive}.
\begin{defn}\label{Definition 2.2}
(Group actions on graphs)
A graph $G=(V, E)$ consists of a set $V$ of vertices and a set $E\subseteq [V]^{2}$ of edges, i.e., $E$ is a subset of the set of all two-element subsets of $V$.
An \emph{action of a group $\Gamma$ on a graph $G=(V,E)$} is a group action of $\Gamma$ on the set $V$ such that, for all $g\in\Gamma$ and all $u,v\in V$,
$
\{u,v\}\in E$ if and only if $\{g u,\, g v\}\in E.
$
For $g\in\Gamma$ and $x\in V$, we denote the action of $g$ on $x$ by $g x$.
This action induces an action of $\Gamma$ on the set of edges $E$ defined by
$
g\{x,y\}=\{g x,\, g y\}$ for all $\{x,y\}\in E.
$
These actions combine to induce a piecewise action of $\Gamma$ on the disjoint union $V \cup E$, defined by $gy = gx$ if $y = x \in V$, and $gy = g\{u, v\}$ if $y = \{u, v\} \in E$.

The {\em action of $\Gamma$ on $G$ is said to be faithful} if the action of $\Gamma$ on $V$ is faithful.
\end{defn}
\begin{defn}\label{Definition 2.3}
Let a group $\Gamma$ act on a graph $G=(V,E)$. For any $Y \subseteq V \cup E$, 
the \emph{stabilizer of $Y$ with respect to the action of $\Gamma$} is defined as
$
\mathrm{Stab}_\Gamma(Y)
=
\{g \in \Gamma : g y = y \text{ for all } y \in Y\}.
$
A subset
$Y\subseteq V\cup E$ is called \emph{$\Gamma$-invariant} if
$g Y=Y$ for every $g\in\Gamma$.
Let $Y\subseteq V\cup E$ be a $\Gamma$-invariant subset and let
$r$ be a positive integer. An \emph{$r$-labeling} of $Y$ is a map
$\ell:Y\to\{1,...,r\}$.
An element $g\in\Gamma$ is said to \emph{preserve} the labeling
$\ell$ if
$\ell(g y)=\ell(y)$ for all $y\in Y$.
The labeling $\ell$ is called \emph{distinguishing (with respect
to the action of $\Gamma$)} if
$
\{g\in\Gamma:\ell(g y)=\ell(y)
\text{ for all } y\in Y\}
=\operatorname{Stab}_\Gamma(Y).
$
\end{defn}

\begin{defn}\label{Definition 2.4}{(Tymoczko \cite{Tym2004})}
The \emph{distinguishing number} $D_\Gamma(X)$ of the action of $\Gamma$ on $X$ is
\[
D_\Gamma(X)
=
\min\{r : \text{there exists a distinguishing $r$-labeling of } X\}.
\]
\end{defn}

Inspired by Definition \ref{Definition 2.4}, we introduce the following definition.

\begin{defn}\label{Definition 2.5}
If $\Gamma$ acts on a graph $G=(V,E)$,  
the elements of $\Gamma$ map vertices to vertices and edges to edges, and so both $E$ and $V\cup E$ are $\Gamma$-invariant. The \emph{distinguishing index} $D'_\Gamma(E)$ of the induced action of $\Gamma$ on $E$ is
    \[
    D'_\Gamma(E)
    =
    \min\{r : \text{there exists a distinguishing $r$-labeling of } E\}
    \]
    and the \emph{total distinguishing number} $D''_\Gamma(V\cup E)$ of the induced action of $\Gamma$ on $V\cup E$ is
    \[
    D''_\Gamma(V\cup E)
    =
    \min\{r : \text{there exists a distinguishing $r$-labeling of } V\cup E\}.
    \]

\end{defn}
\begin{defn}\label{Definition 2.6}
A group $\Gamma$ is called a \emph{semidirect product} of $N$ by $Q$, 
denoted by $\Gamma = N \rtimes Q$, if $\Gamma$ contains subgroups $N$ and $Q$ such that:
(1).  $N \trianglelefteq \Gamma$ (that is, $N$ is a normal subgroup of $\Gamma$), 
(2). $NQ = \Gamma$, and
(3). $N \cap Q = \{e\}$.
Let $\Gamma$ be a group and $S_m$ be the symmetric group on 
$\{1,...,m\}$. The \emph{wreath product} of $\Gamma$ by $S_m$, denoted $\Gamma \wr S_m$, 
is defined as the semidirect product
$
\Gamma \wr S_m = \Gamma^m \rtimes S_m,
$
where $\Gamma^m$ denotes the direct product of $m$ copies of $\Gamma$, and $S_m$ acts on $\Gamma^m$ by permuting coordinates.
An element of $\Gamma \wr S_m$ has the form
$
(g_1,g_2,...,g_m;\sigma),
$
where $g_i \in \Gamma$ for $1 \le i \le m$ and $\sigma \in S_m$. The group operation is 
given by
$
(g_1,...,g_m;\sigma)(h_1,...,h_m;\tau)
=
(g_1h_{\sigma^{-1}(1)},...,g_mh_{\sigma^{-1}(m)};\sigma\tau).
$
Let $\Gamma$ be a group acting on a set $X$. The wreath product $\Gamma \wr S_m$ acts naturally on the Cartesian power $X^m$ by
$
(g_1,...,g_m;\sigma)(x_1,...,x_m)
=
(g_1x_{\sigma^{-1}(1)},...,g_mx_{\sigma^{-1}(m)}),
$
for all $(x_1,...,x_m)\in X^m$.
\end{defn}

\begin{defn}\label{Definition 2.7}
Let $\Gamma$ be a group and $S \subseteq \Gamma \setminus \{e\}$ 
be inverse-closed i.e., $S = S^{-1}$, where 
$S^{-1}:= \{s^{-1}: s \in S\}$.
The \emph{undirected Cayley graph} Cay$(\Gamma, S)$ is the graph with a set of vertices $\Gamma$, and the vertices $u$ and $v$ are adjacent in Cay$(\Gamma, S)$ if and only if $uv^{-1} \in S$. It is known that Cay$(\Gamma, S)$ is connected if and only if $S$ is a generating set of $\Gamma$.
\end{defn}

\begin{defn}\label{Definition 2.8}
The \emph{right regular representation} of a group $\Gamma$, denoted by $R(\Gamma)$,
is the permutation group 
$\{\rho_g : \Gamma \to \Gamma \mid \rho_g(x)=xg 
\text{ for all }
x\in \Gamma,\; g\in \Gamma\}$ in $\mathrm{Sym}(\Gamma)$.
It is known that $R(\Gamma)$ is a subgroup of $\mathrm{Aut}(\mathrm{Cay}(\Gamma,S))$ and if 
$R(\Gamma)=\mathrm{Aut}(\mathrm{Cay}(\Gamma,S))$, then 
$\mathrm{Cay}(\Gamma,S)$ is a {\em graphical regular representation} of $\Gamma$.
\end{defn}

\begin{defn}\label{Definition 2.9}
Given a group $\Gamma$ and a subset $S \subseteq \Gamma$, we define 
$\mathrm{Aut}(\Gamma,S) = \{\alpha \in \mathrm{Aut}(\Gamma) \mid \alpha(S) = S\}$.
A Cayley graph $G=$ Cay$(\Gamma,S)$ is {\em normal} if $R(\Gamma)\trianglelefteq \mathrm{Aut}(G)$. The graph $G$ is normal if and only if $\mathrm{Aut}(G)=R(\Gamma) \rtimes \mathrm{Aut}(\Gamma,S)$.
\end{defn}

\begin{defn}[Cartesian power of a Graph]\label{Definition 2.10}
Let $G = (V, E)$ be a graph and $k \ge 1$ be an integer. The \emph{Cartesian power} of $G$, denoted by $G^k$, is the graph defined by the vertex set $V(G^k) = V^k$ and the edge set $E(G^k) = E_k$, where:
\begin{enumerate}
    \item 
    $
    V^k = \prod_{i=1}^k V = \{(u_1, u_2, ..., u_k) \mid u_i \in V \text{ for all } 1 \le i \le k\}.
    $
    \item Two vertices $\mathbf{u} = (u_1, ..., u_k)$ and $\mathbf{v} = (v_1, ..., v_k)$ in $V^k$ are adjacent in $G^k$ if and only if there exists a unique index $\ell \in \{1, ..., k\}$ such that $\{u_\ell, v_\ell\} \in E$ and $u_i = v_i$ for all $i \neq \ell$.
\end{enumerate}
\end{defn}

\section{Cayley graphs}

\begin{fact}\label{Fact 3.1}
{\em The following hold:
\begin{enumerate}
    \item (Lehner--Pil\'{s}niak--Stawiski; \cite[Theorem 3]{LPS2020}) 
    Let $G$ be a connected $\Delta$-regular graph for $\Delta \in \mathbb{N} \setminus \{1\}$. Then $D'(G) \leq 3$.
    
    \item (Alikhani--Soltani; \cite[Theorem 2.3, 3.2(iii)]{AS2021})
    If $G = Cay(H, S)$ is a Cayley graph of order $|H| \geq 2$ and $G$ is a graphical regular representation of $H$ then $D'(G)=D(G)= 2$.
   
    \item (Alikhani--Soltani; \cite[Theorem 2.2]{AS2021})
    If $G = \mathrm{Cay}(H, S)$ is normal, then $D(G) \leq 1+ D_{\mathrm{Aut(H,S)}}(S)$.   
\end{enumerate}
}
\end{fact}

Similar to Fact \ref{Fact 3.1}(3), we obtain upper bounds of $D''(G)$ and $D'(G)$ with respect to $D''_{\mathrm{Aut(H,S)}}(S\cup E_{e}\cup \{e\})$ and $D'_{\mathrm{Aut(H,S)}}(E_{e})$ respectively where $E_{e}=\{\{e,s\}:s\in S\}$ and $e$ is the identity of $H$.

\begin{thm}\label{Theorem 3.2}
{\em Let $G=\mathrm{Cay}(H,S)$ be a Cayley graph on a
group $H$ with identity $e$, where $|H|\geq 2$ and $|S|\geq 2$. The following hold:
\begin{enumerate}
    \item If $G$ is connected then $2\leq D''(G)\leq D'(G)\leq 3$;
    \item If $G$ is a graphical regular representation of $H$, then
    $D''(G)=2$.
    \item If $G$ is normal (not necessarily connected) such that
    \begin{enumerate}
        \item[(i)] $\operatorname{Stab}_{R(H)}^{\mathrm{set}}(S)
        =\{R_a\in R(H):Sa=S\}
        =\{\operatorname{id}_{R(H)}\}$ where $Sa=\{sa:s\in S\}$, and
        \item[(ii)] $\operatorname{Stab}_{\operatorname{Aut}(H,S)}(S)
        =\{\operatorname{id}_{\operatorname{Aut}(H,S)}\}$,
    \end{enumerate}
    then
    $
    D''(G)\leq
    1+D''_{\operatorname{Aut}(H,S)}
    (S\cup E_e\cup\{e\}),
    $
    where $E_e=\{\{e,s\}:s\in S\}$.
\end{enumerate}
}
\end{thm}

\begin{proof}
(1). It is clear that $D''(G) \leq D'(G)$. Since $R(H)$ is a subgroup of $\mathrm{Aut}(G)$, 
we have
$|\operatorname{Aut}(G)|\geq |H|\geq 2$.
Thus,
$\mathrm{Aut}(G)$ has a non-trivial automorphism, and so $D''(G) \geq 2$. On the other hand, $G$ is a connected $|S|$-regular graph.
Thus, $D'(G) \leq 3$ by Fact \ref{Fact 3.1}(1).     

(2). Since $\mathrm{Aut}(G)=R(H)$ and $|H|\geq 2$, so
$D''(G)\ge 2$. Since $D''(G)\leq \mathrm{min}\{D(G),D'(G)\}=2$ by Fact \ref{Fact 3.1}(2), we obtain $D''(G)=2$.

(3). 
Let $\operatorname{Aut}(G)$ act on $G$.
Since every $g\in\operatorname{Aut}(H,S)$ fixes $e$ and maps $S$ onto itself, it also maps $E_e$ onto itself. Thus, $S\cup E_e\cup\{e\}$ is $\operatorname{Aut}(H,S)$-invariant under the induced action of $\operatorname{Aut}(H,S)$ on $V\cup E$.
Let $t= D''_{\mathrm{Aut}(H,S)}(S\cup E_{e}\cup \{e\})$ and $\ell:S\cup E_{e}\cup \{e\}\rightarrow \{1,...,t\}$ be a total distinguishing labeling 
for the action of $\mathrm{Aut}(H,S)$ on $S\cup E_e\cup \{e\}$ (i.e. if $g \in \mathrm{Aut}(H,S)$ preserves $\ell$, then $g\in \mathrm{Stab}_{\mathrm{Aut}(H,S)}(S\cup E_{e}\cup \{e\})$).
We extend $\ell$ to a total labeling $\mathcal L$ of the graph $G$
using the labels
$0,1,...,t$ as follows:
\begin{itemize}
  \item Let $S=\{s_1,...,s_k\}$. For each $s\in S\cup \{e\}$, define $\mathcal L(s):=\ell(s)$. 
  \item For each edge $\{e,s\}\in E_e$, define $\mathcal L(\{e,s\}):=\ell(\{e,s\})$.
  \item Assign a new label $0$ to every remaining vertices and edges of $G$.
\end{itemize}

It suffices to show that $\mathcal{L}$ is distinguishing. 
Let $\varphi\in\mathrm{Aut}(G)$ be any $\mathcal L$-preserving automorphism. 

\begin{claim}\label{Claim 3.3}
    {\em $\varphi(S) = S$.}
\end{claim}

\begin{proof}
Since the labels $1,...,t$ occur only on elements of $S\cup E_e\cup\{e\}$,  $\varphi(S\cup E_e\cup\{e\}) = S\cup E_e\cup\{e\}$.
Since $\varphi$ sends vertices to vertices, $\varphi(S \cup \{e\}) = S \cup \{e\}$.
Now, $|S|\ge 2$. So, the identity element $e$ of $H$ is the unique vertex incident to every edge in $E_e$, and hence it is fixed by $\varphi$. Thus, $\varphi(S) = S$.    
\end{proof}

Since $G$ is normal, $\mathrm{Aut}(G)=R(H)\rtimes \mathrm{Aut}(H,S)$.
Thus, there exist $g \in \mathrm{Aut}(H, S)$ and $a \in H$ such that $\varphi(x) = g(x)a$ for all $x \in H$, i.e., $\varphi = R_a \circ g$, where $R_a\in R(H)$ is defined by $R_a(x)=xa$ for all $x\in H$.

\begin{claim}\label{Claim 3.4}
{\em $R_a = \mathrm{id}_{R(H)}$ and $g=\mathrm{id}_{\mathrm{Aut}(H,S)}$.} 
\end{claim}

\begin{proof}
Since $g \in \mathrm{Aut}(H, S)$, we have $g(S) = S$. 
Thus, $\varphi(S)  = \{\varphi(s) : s \in S\} = \{g(s)a : s \in S\} = \{sa : s \in S\} = Sa$. 
By Claim \ref{Claim 3.3}, $Sa = S$. So, $R_a \in \mathrm{Stab}_{R(H)}^{\mathrm{set}}(S)$. By the assumption $\mathrm{Stab}_{R(H)}^{\mathrm{set}}(S) = \{\mathrm{id}_{R(H)}\}$, it follows that $R_a = \mathrm{id}_{R(H)}$.
Since $g$ preserves $\ell=\mathcal{L}\!\restriction (S\cup E_e\cup\{e\})$ and $g\in \mathrm{Aut}(H,S)$, $g\in \mathrm{Stab}_{\mathrm{Aut}(H,S)}(S\cup E_{e}\cup\{e\})$.
Since $S \subseteq S \cup E_e \cup \{e\}$, we obtain
$g\in \mathrm{Stab}_{\mathrm{Aut}(H,S)}(S)$.
By the assumption 
$\mathrm{Stab}_{\mathrm{Aut}(H,S)}(S)
= \{\mathrm{id}_{\mathrm{Aut}(H,S)}\}$, we obtain $g = \mathrm{id}_{\mathrm{Aut}(H,S)}$.
\end{proof}

By Claim \ref{Claim 3.4}, 
$\varphi = \mathrm{id}_{R(H)} \circ
\mathrm{id}_{\mathrm{Aut}(H,S)} = \mathrm{id}_{G}$.\footnote{
In the semidirect product decomposition $\mathrm{Aut}(G)=R(H)\rtimes \mathrm{Aut}(H,S)$, each element of 
$\mathrm{Aut}(G)$ is written uniquely as $R_a \circ g$ 
with $a\in H$ and $g\in \mathrm{Aut}(H,S)$ and the identity element 
of $\mathrm{Aut}(G)$ is
$\mathrm{id}_{R(H)} \circ \mathrm{id}_{\mathrm{Aut}(H,S)}$.}
Thus, $\mathcal{L}$ is
a distinguishing labeling.
\end{proof}

Similar to Theorem \ref{Theorem 3.2}(3), we can prove the following.
\begin{thm}\label{Theorem 3.5}
{\em 
Let $G=\mathrm{Cay}(H,S)$ be a normal Cayley graph with $|S| \geq 2$ where $\mathrm{Stab}_{R(H)}^{\mathrm{set}}(S)=\{\mathrm{id}_{R(H)}\}$ and 
$\mathrm{Stab}_{\mathrm{Aut}(H,S)}(S)=\{\mathrm{id}_{\mathrm{Aut}(H,S)}\}$.
Then $D'(G)\le 1 + D'_{\mathrm{Aut}(H,S)}(E_e)$ where
$E_e=\{\{e,s\}:s\in S\}$ and $e$ is the identity of $H$.
}
\end{thm}
\section{Elementary Group actions}

\begin{prop}\label{Proposition 4.1}
{\em For any action of a group $\Gamma$ on a graph $G = (V, E)$, we have:
\begin{enumerate}
    \item $\mathrm{Stab}_\Gamma(V) = \mathrm{Stab}_\Gamma(V \cup E)$,
    \item $D''_\Gamma(V \cup E) \leq D_\Gamma(V)$, and
    \item $D''_\Gamma(V \cup E) \leq D'_\Gamma(E)$ if 
    $G$ is a connected graph of order at least $3$.
\end{enumerate}
}
\end{prop}

\begin{proof}
(1). 
If $g\in \mathrm{Stab}_\Gamma(V\cup E)$, then $g\in \mathrm{Stab}_\Gamma(V)$ as $V\subseteq V\cup E$.
Let $g\in \mathrm{Stab}_\Gamma(V)$.
Pick any edge $\{u,v\}\in E$. 
Since $u, v \in V$ and $g \in \mathrm{Stab}_{\Gamma}(V)$, we have $g u = u$ and $g v = v$.
Under
the induced action of $\Gamma$ on $E$, we have
$
g\{u,v\}
=
\{g u,g v\}
=
\{u,v\}
$.
Hence $g$ fixes every edge of $E$, and so
$
g\in \mathrm{Stab}_\Gamma(V\cup E).
$

(2). Let $d = D_\Gamma(V)$. Then, there exists a distinguishing $d$-labeling $\phi : V \to \{1, ..., d\}$. 
We extend $\phi$ to a total labeling $\tau : V \cup E \to \{1, ..., d\}$ by defining $\tau(y)=\phi(y)$ if $y\in V$ and $\tau(y)=1$ if $y\in E$.
If $g \in \Gamma$ preserves $\tau$, then $g$ preserves $\phi=\tau\restriction V$. 
Since $\phi$ is distinguishing, by (1), we have $g\in \mathrm{Stab}_\Gamma(V) = \mathrm{Stab}_\Gamma(V \cup E)$. Thus, $\tau$ is distinguishing, proving $D''_\Gamma(V \cup E) \leq D_\Gamma(V)$. 

(3). In view of the arguments in (2), it suffices to prove $\operatorname{Stab}_\Gamma(E)
=\operatorname{Stab}_\Gamma(V\cup E)$. Since $E\subseteq V\cup E$,
$\operatorname{Stab}_\Gamma(V\cup E)
\subseteq \operatorname{Stab}_\Gamma(E)$.
Let $g\in\operatorname{Stab}_\Gamma(E)$.
Pick any vertex $x\in V$. Since $G$ is connected of order at least $3$,
$\deg(x)\geq 1$ where $\deg(x)$ is the degree of $x$.

\textbf{Case 1:} $\deg(x)\geq 2$. Let $\{x, y\}$ and $\{x, z\}$ be two distinct edges incident to $x$. Since $g\in Stab_{\Gamma}(E)$, $g\{x, y\}=\{gx,gy\}=\{x, y\}$ and $g\{x, z\}=\{gx,gz\}=\{x, z\}$. Thus, $gx\in \{x,y\}\cap \{x,z\}=\{x\}$, so $gx=x$.    

\textbf{Case 2:} $\deg(x)=1$.
Let $y$ be the unique neighbor of $x$. Since $G$ is connected and $|V|\geq 3$, we have $\deg(y)\geq 2$. By the arguments in Case 1, $gy=y$. Since $g\in Stab_{\Gamma}(E)$, $\{gx,gy\}=\{x,y\}$. Thus, $gx=x$.

Thus, $g$ fixes every vertex, and hence
$g\in\operatorname{Stab}_\Gamma(V\cup E)$.
Therefore,
$\operatorname{Stab}_\Gamma(E)
=\operatorname{Stab}_\Gamma(V\cup E)$.
\end{proof}

The following proposition follows from the definitions.

\begin{prop}\label{Proposition 4.2}
{\em Let $\Gamma$ be a group such that $|\Gamma|>1$. Then $\Gamma$ acts on the nontrivial graph $G=(V,E)$ by fixing
every vertex if and only if
$D''_{\Gamma}(V\cup E)=1$.}
\end{prop}

\subsection{Base of size $n$}

\begin{lem}(\cite[Lemma 13]{Ali2026})\label{Lemma 4.3}
{\em If a group $\Gamma$ acts on a set $V$ with a base of size $n$, then $D_{\Gamma}(V)\le n+1$.} 
\end{lem}
Suppose the vertex set of $K_{1,n}$ is given by $V = \{v_0, ..., v_n\}$ with central vertex $v_0$, and the edge set is $E = \{e_1, ..., e_n\}$ with $e_i = \{v_0, v_i\}$. 
Let $\mathrm{id}_{S_n}$ and $\mathrm{id}_{A_n}$ be the identities of $S_n$ and $A_n$ respectively.

\begin{thm}\label{Theorem 4.4}
{\em Let a group $\Gamma$ act on a graph $G=(V,E)$.
The following hold:
\begin{enumerate}
    \item If the action of $\Gamma$ on $V$ has a base of size $n$, then
    $D''_{\Gamma}(V\cup E)\le n+1$.
    \item If the induced action of $\Gamma$ on $E$ has a base of size $m$, then
    $D'_{\Gamma}(E)\le m+1$.

    \item The upper bound in (2) is sharp for the following actions:
\begin{enumerate}[(i)]
    \item The action of $S_n$ on $K_{1,n}$ for $n\ge 3$, defined by
$\sigma(v_0)=v_0$, $\sigma(v_i)=v_{\sigma(i)}$ for
$1\le i\le n$.
    \item The action of $A_n$ on $K_{1,n}$ for $n\ge 5$, defined by
$\sigma(v_0)=v_0$, $\sigma(v_i)=v_{\sigma(i)}$ for
$1\le i\le n$.
\end{enumerate}
\end{enumerate}
}
\end{thm}

\begin{proof}
(1). Let $Y=V\cup E$. Let $B=\{v_1,...,v_n\}\subseteq V$ be a base for the action of $\Gamma$ on $V$, which means 
$
\bigcap_{i=1}^n \mathrm{Stab}_{\Gamma}(v_i)=\{e\}.
$
Define $\ell:Y\to\{1,...,n+1\}$ by
$
\ell(v_i)=i$ for all $1\le i\le n$,
and assign a new label $n+1$ to every element of $Y \setminus B$.
Suppose $g\in \Gamma$ preserves $\ell$. Since $v_i$ is the unique element in $Y$ with label $i$, we have $g(v_i)=v_i$ for all $1\le i\le n$. 
Hence, 
$
g\in\bigcap_{i=1}^n \mathrm{Stab}_{\Gamma}(v_i)=\{e\},$
which implies $g=e\in\operatorname{Stab}_\Gamma(Y)$.
Thus $\ell$ is a total distinguishing labeling, proving
$D''_{\Gamma}(Y)\le n+1$. 

(2). Let $B=\{e_1,...,e_m\}$ be a base for the action on $E$.
Define a labeling $\ell:E\to\{1,...,m+1\}$ by assigning
label $i$ to $e_i$ for each $1\le i\le m$, and label $m+1$
to every edge in $E\setminus B$.
The rest of the proof follows analogously to the proof of (1).

(3).
Since every element of $S_n$ (respectively, $A_n$) fixes
$v_0$ and permutes the leaves, the induced action of $S_n$
(respectively, $A_n$) on the edge set $E$ is their natural
action on the set $\{1,...,n\}$.

$(i)$. 
We prove that the base size for the action of $S_n$ on $E$ is
$n-1$. If $B\subseteq E$ has size $n-1$, say
$B=\{e_1,...,e_{n-1}\}$, then any
$g\in\bigcap_{b\in B}\operatorname{Stab}_{S_n}(b)$ fixes
$e_1,...,e_{n-1}$. So, the remaining edge $e_n$ is
also fixed, and  $g=\mathrm{id}_{S_n}$, confirming that
$B$ is a base. However, no $B\subseteq E$ of size $n-2$ is
a base, since the transposition $(i\,j)\in S_n$, where
$e_i,e_j\notin B$, interchanges $e_i$ and $e_j$ under the
induced action on $E$ and fixes every edge in $B$.
However, $D'_{S_n}(E) = D'_{\operatorname{Aut}(K_{1,n})}(E) = D'(K_{1,n}) = n$.


$(ii)$.
We show that the base size for the action of $A_n$ on $E$ is
$n-2$. Let $B\subseteq E$ be of size $n-2$, say
$B=\{e_1,...,e_{n-2}\}$.
Any $g\in\bigcap_{b\in B}\operatorname{Stab}_{A_n}(b)$
fixes $e_1,...,e_{n-2}$. Hence, either
$g=\mathrm{id}_{A_n}$ or $g$ is the transposition $(n-1\ n)$.
Since $(n-1\ n)\notin A_n$, we get $g=\mathrm{id}_{A_n}$. Thus, $B$ is a base.
However, any subset $B'$ of $E$ of size at most
$n-3$ leaves at least three edges outside $B'$. A $3$-cycle
on the indices corresponding to three such edges is a
nonidentity element of $A_n$ fixing every edge in $B'$. Hence, $B'$ is not a base.

Now, we show that $D'_{A_n}(E)=n-1$. By (2),
$D'_{A_n}(E)\le n-2+1=n-1$. Since $A_n$ acts faithfully
on $E$, $\operatorname{Stab}_{A_n}(E)=\{\mathrm{id}_{A_n}\}$.
Suppose we use at most $n-2$ labels on $E$. Then any
such labeling of the $n$ edges contains either a label class
of size at least $3$, or at least two label classes of size
$2$.

\begin{itemize}
    \item[] \textit{Case 1:} A label class contains at least three distinct
    edges, say $e_i,e_j,e_k$. The $3$-cycle $g=(i\,j\,k)=(i\,k)(i\,j)\in A_n\setminus\{\mathrm{id}_{A_n}\}$
    permutes $e_i,e_j$ and $e_k$, and fixes all other
    edges under the induced action on $E$. Since $e_i,e_j$ and $e_k$ 
    have the same label, $g$ preserves the labeling.
\vspace{2mm}

    \item[] \textit{Case 2:} 
    At least two label classes have size two, say $\{e_i, e_j\}$ and $\{e_k, e_l\}$ with $i, j, k, l$ distinct.
    Then $g=(i\,j)(k\,l)\in A_n\setminus\{\mathrm{id}_{A_n}\}$ interchanges the edges
    within each label class and fixes all other edges under the
    induced action on $E$. So, $g$ preserves the labeling.
\end{itemize}

In both the cases, we have found a nonidentity element of $A_n$
that preserves the labeling. Hence, the labeling is not
distinguishing. Thus, $D'_{A_n}(E)\ge n-1$, and therefore
$D'_{A_n}(E)=n-1$.
\end{proof}
\subsection{Distinguishing invariants under the Action of Subgroups}

\begin{prop}\label{Proposition 4.5}
{\em The following holds:
\begin{enumerate}
    \item Suppose $\Gamma$ acts faithfully on a graph $G=(V,E)$ and $H\leq \Gamma$. Consider the action
    of $H$ on $G$ obtained by restricting the action of $\Gamma$. Then $D'_H(E)\leq D'_\Gamma(E)$ and $D''_{H}(V\cup E)\leq D''_{\Gamma}(V\cup E)$.
    \item For any integer $n\geq 2$,
    consider the action
    of the rotational subgroup $\langle\rho\rangle$ of the dihedral group $D_{n}=\mathrm{Aut}(C_n)$ on $C_n$ obtained by restricting the action of $D_{n}$. Then $D'_{\langle\rho\rangle}(E(C_n))=2$.
    
    \item For any integer $n\geq 2$, consider the action of $A_n\leq S_n=\mathrm{Aut}(K_n)$ on $K_n$ obtained by restricting the action of $S_{n}$. Then $D'_{A_n}(E(K_n))=2$.
    
    \item There exist infinitely many graphs $G$ with
$\Gamma\leq\operatorname{Aut}(G)$ where
$D'_{\Gamma}(E(G))< D'(G)$.
\end{enumerate}
}
\end{prop}

\begin{proof}
(1). Assume that $\ell$ is a distinguishing edge labeling of $G$
with $D'_\Gamma(E)$ labels. 
Since $\Gamma$ has no nonidentity element preserving $\ell$, neither does $H \le \Gamma$.
So $D'_{H}(E)$ colors suffice to
distinguish the action of $H$ on $G$. Thus, $D'_H(E)\leq D'_\Gamma(E)$. A similar argument shows $D''_{H}(V\cup E)\leq D''_{\Gamma}(V\cup E)$.

(2). Since $|\langle\rho\rangle|>1$, we have
$D'_{\langle\rho\rangle}(E(C_n))\geq2$.
Assign one label to an edge
of $C_n$ and another label to all the remaining edges.
Every nonidentity rotation moves the uniquely labeled edge.
Thus, only the identity of $\langle\rho\rangle$ preserves this labeling, and hence
$D'_{\langle\rho\rangle}(E(C_n))=2$.

(3). Since $|A_n|>1$, we have $D'_{A_n}(E(K_n))\geq2$.
Let $V(K_n)=\{1,...,n\}$.
If $n=3$, we assign label $a$ to the edges $\{1,3\}$ and $\{2,3\}$, and  label $b$ to the edge
$\{1,2\}$. The only nonidentity permutation preserving this labeling
is $(1\,2)$ and $(1\,2)\not\in A_n$. Thus, the labeling
is distinguishing for the action of $A_n$.
Now let $n\geq4$. Assign label $a$ to the edges
$\{1,2\},\{1,3\},\{2,3\},\{3,4\},\{4,5\},...,\{(n-1),n\}$, and label $b$ to all the
remaining edges of $K_n$.
The subgraph $H$ induced by the edges with label $a$ contains the triangle $123$ with a path attached at vertex $3$. Let $g\in A_n$ preserves this edge labeling
of $K_n$. Then $g$ preserves the set of edges of $H$.
However, the only nonidentity element of $S_n$ that preserves
the edge set of $H$ is the transposition $(1\,2)$.
Since $(1\,2)\notin A_n$, $g=\mathrm{id}_{A_n}$. So,
$D'_{A_n}(E(K_n))= 2$.

(4). Fix $n\in\{3,4,5\}$. 

\textbf{First example.} If the rotational subgroup $\langle\rho\rangle$ of the dihedral group $D_{n}=\mathrm{Aut}(C_n)$ acts on $C_n$ as in (2), then $D'_{\langle\rho\rangle}(E(C_n))=2<D'(C_n)=3$.

\textbf{Second example.} Consider the action of $A_{n}\leq S_n=\mathrm{Aut}(K_n)$ on $K_n$ as in (3). Then, $D'_{A_n}(E(K_n))= 2<D'(K_n)=3$.

\textbf{An infinite family.} Consider the action of $A_k\leq S_k=\mathrm{Aut}(K_{1,k})$ on $K_{1,k}$ for any integer $k\ge 5$, defined by
$\sigma(v_0)=v_0$, $\sigma(v_i)=v_{\sigma(i)}$ for
$1\le i\le k$. As established in the proof of Theorem \ref{Theorem 4.4} $(ii)$, $D'_{A_k}(E(K_{1,k}))= k-1<D'(K_{1,k})=k$.
Since the graphs $K_{1,k}$, for $k \geq 5$, are pairwise non-isomorphic, this establishes the assertion for infinitely many graphs.
\end{proof}
\begin{prop}\label{Proposition 4.6}
{\em The following holds:
\begin{enumerate}
    \item For any integer $k\geq 5$, $D''_{A_k}(V(K_{1,k})\cup E(K_{1,k}))=\lceil \sqrt{k-1} \rceil$ if $A_k$ acts faithfully on the star $K_{1,k}$ by
$\sigma(v_0)=v_0$, $\sigma(v_i)=v_{\sigma(i)}$ for
$1\le i\le k$.
    
    \item There exist infinitely many graphs $G$ with
$\Gamma\leq\operatorname{Aut}(G)$ where
$D''_{\Gamma}(V(G)\cup E(G))< D''(G)$.
\end{enumerate}
}
\end{prop}

\begin{proof}
(1). Let $q=\lceil \sqrt{k-1} \rceil$ and $X=\{(v_i,e_i): v_i \text{ is a leaf}, 1\leq i \leq k\}$ where $e_i=v_0v_i$. Since $q\geq2$, we have at least two
available labels. We color the central vertex $v_0$ with $2$.  
Assign $(1,1)$ to $(v_1,e_1)$ and $(v_2,e_2)$ and label each $k-2$ pair of $X\backslash \{(v_1,e_1),(v_2,e_2)\}$ with a distinct pair of labels different from $(1,1)$. 
This is possible since $q^2\geq k-1$. Suppose $g\in A_k$ preserves this total labeling. Thus, $g=\mathrm{id}_{A_k}$ or $g=(1\ 2)$. Since $(1\ 2)\not\in A_k$, $g=\mathrm{id}_{A_k}$. Thus, $D''_{A_k}(V(K_{1,k})\cup E(K_{1,k}))\leq q$. 
Consider any total labeling $\ell$ of $K_{1,k}$ with at most $q-1$ labels. Then the number of label pairs $(q-1)^2\leq n-2$. Thus, at least one label class contains at least $3$ elements of $X$, say $(v_i,e_i), (v_j,e_j), (v_k,e_k)$ or at least two label classes contain $2$ elements of $X$. 
In the first case, the permutation $(i\ j\ k)\in A_n\backslash \{\mathrm{id}_{A_n}\}$ preserves $\ell$. In the second case, suppose $(v_i,e_i), (v_j,e_j)$ have one label and  $(v_k,e_k), (v_l,e_l)$ have another label. Then, the permutation $(i\ j) (k\ l)\in A_n\backslash \{\mathrm{id}_{A_n}\}$ preserves $\ell$.
Thus, $\ell$ is not distinguishing. So, $D''_{A_k}(V(K_{1,k})\cup E(K_{1,k}))\geq q$.

(2). Consider $k=t^2+1$ for any  $t\geq 2$. 
Let $A_k\leq S_k=\mathrm{Aut}(K_{1,k})$ acts faithfully on $K_{1,k}$. By (1), $
D''_{A_k}(V(K_{1,k})\cup E(K_{1,k}))=\lceil \sqrt{k-1} \rceil=t$, whereas any distinguishing total labeling of
$K_{1,k}$ under its full automorphism group $S_k$ must
assign distinct pairs of labels to all $k$ elements of $X$.
So,
$D''(K_{1,k})=\lceil\sqrt{k}\rceil=t+1$.
Since $K_{1,k}$, for $k=t^2+1$, $t\geq2$,
are pairwise non-isomorphic, the assertion follows.
\end{proof}
\begin{thm}\label{Theorem 4.7}
{\em For every $n\geq 3$,
$
D'_{A_3\times S_n}\bigl(E(K_{3,n})\bigr)
=\left\lceil (n+1)^{1/3}\right\rceil$.
}
\end{thm}

\begin{proof}
Let $X=\{x_1,x_2,x_3\}$ and
$Y=\{y_1,...,y_n\}$ be the bipartitions of $K_{3,n}$.
The group $A_3\times S_n$ acts faithfully on $K_{3,n}$, where
$A_3$ permutes the vertices of $X$ and $S_n$ permutes the
vertices of $Y$.

\begin{claim}\label{Claim 4.8}
$D'_{A_3\times S_n}\bigl(E(K_{3,n})\bigr)
\geq\left\lceil (n+1)^{1/3}\right\rceil$.
\end{claim}

\begin{proof}
Suppose that $f$ is a
distinguishing edge labeling of $K_{3,n}$ using $q$ labels, and let
$\Omega=\{0,...,q-1\}^3$. For each $j\in\{1,...,n\}$,
define the tuple
$
c_j=\bigl(f(x_1y_j),f(x_2y_j),f(x_3y_j)\bigr)\in\Omega.
$

\textbf{Subclaim.} {\em The tuples $c_1,...,c_n$ must be pairwise distinct.}

\textit{Proof.} 
If $c_i=c_j$ for some $i\ne j$, then the non identity element
$(\mathrm{id}_{A_3},(i\,j))\in A_3\times S_n$ preserves $f$ and
interchanges $y_i$ and $y_j$,
contrary to the assumption that $f$ is distinguishing. 

We show that the set
$S=\{c_j:1\leq j\leq n\}$
cannot be $\Omega$. 
Suppose $S=\Omega$. 
The action of $(123)\in A_3$ on $X$ induces the permutation
$\rho:\Omega\rightarrow\Omega$ where $\rho(a,b,c)=(b,c,a)$.
Since $S=\Omega$, we have $\rho(S)=S$. Since $c_1,...,c_n$ are pairwise
distinct and comprise all elements of $S=\Omega$, for each
$j\in\{1,...,n\}$ there exists a unique index
$\tau(j)\in\{1,...,n\}$ such that
$c_{\tau(j)}=\rho^{-1}(c_j)$.
Since $\rho^{-1}$ is a bijection of $\Omega$, the map
$\tau$ is a permutation of $\{1,...,n\}$. Now
$c_j=\rho(c_{\tau(j)})$
for every $j\in \{1,...,n\}$. Writing
$c_j=(a_j,b_j,d_j)$, this yields
$
(a_j,b_j,d_j)
=(b_{\tau(j)},d_{\tau(j)},a_{\tau(j)}).
$
Therefore, 
\[f(x_1y_j)=f(x_2y_{\tau(j)}),
f(x_2y_j)=f(x_3y_{\tau(j)}), 
\text{ and }
f(x_3y_j)=f(x_1y_{\tau(j)})
\]
for every $j$.
Hence the nonidentity element
$((123),\tau)\in A_3\times S_n$ preserves $f$, a
contradiction. Thus $S\ne\Omega$, and consequently
$n\leq q^3-1$, and therefore
$q\geq \left\lceil (n+1)^{1/3}\right\rceil$.
\end{proof}

\begin{claim}\label{Claim 4.9}
$D'_{A_3\times S_n}\bigl(E(K_{3,n})\bigr)
\leq\left\lceil (n+1)^{1/3}\right\rceil$.
\end{claim}

\begin{proof}
Let
$q=\left\lceil (n+1)^{1/3}\right\rceil$ and $\Omega=\{0,...,q-1\}^3$.
Then $q\geq 2$ and $n\leq q^3-1$.
Consider the permutation
\[
\rho:\Omega\rightarrow\Omega,\qquad
\rho(a,b,c)=(b,c,a).
\]
Since $\rho^3$ is the identity, 
the cyclic group
$\langle\rho\rangle=\{\mathrm{id},\rho,\rho^2\}$ acts on $\Omega$.
For each $u\in\Omega$, its orbit under this action is
$\operatorname{Orb}_{\langle\rho\rangle}(u)
=\{g(u):g\in\langle\rho\rangle\}
=\{u,\rho(u),\rho^2(u)\}$.
Consequently, every orbit has size either $1$ or $3$. The orbits
of size $1$ are the singleton sets
$\{(a,a,a)\}$, where $a\in\{0,...,q-1\}$, while every other
orbit has size $3$. Since $q\ge2$, there exists at least one
orbit of size $3$. Choose such an orbit
$O=\{u,\rho(u),\rho^2(u)\}$.

\textbf{Step 1.} We choose an $n$-element subset $S\subseteq\Omega$ such that $\rho(S)\neq S$ and $\rho^2(S)\neq S$.

\begin{itemize}
    \item []Since $|O|=3$, we have
    $|\Omega\setminus O|=q^3-3\geq n-2$.
    Hence, we can choose $n-2$ elements from $\Omega\setminus O$ and define $S$ to be their union with $\{u,\rho(u)\}$.
    Thus, $|S|=n$, $S\cap O=\{u,\rho(u)\}$, and
    $\rho^2(u)\notin S$.
    If $\rho(S)=S$, then, since $\rho(u)\in S$,
    we would have $\rho^2(u)\in S$, a contradiction.
    Similarly, if $\rho^2(S)=S$, then, since $u\in S$,
    we would have $\rho^2(u)\in S$, again a contradiction.
\end{itemize}

\textbf{Step 2.}
We define an edge labeling $f$ of $K_{3,n}$.

\begin{itemize}
    \item []Let
$c_{j}=(a_j,b_j,d_j)\in S$
for all $1\leq j\leq n$ be $n$ distinct elements of $S$.
Define an edge labeling $f$ of $K_{3,n}$
by $f(x_1y_j)=a_j,
f(x_2y_j)=b_j$,
and $f(x_3y_j)=d_j$
for every $j\in\{1,...,n\}$.
\end{itemize}

We will show that $f$ is distinguishing.
Suppose that
$(\sigma,\tau)\in A_3\times S_n$ preserves $f$. 

\textbf{Step 3.}
We prove $\sigma=\mathrm{id}_{A_3}$.

\begin{itemize}
    \item [] We note that
$f(x_i y_j)=f\bigl(\sigma(x_i)y_{\tau(j)}\bigr)
$
for every $i\in\{1,2,3\}$ and every
$j\in\{1,...,n\}$.
Suppose that $\sigma=(123)$. 
Then $\sigma(x_1)=x_2, \sigma(x_2)=x_3,$ and $\sigma(x_3)=x_1$. Thus,
\[
\begin{aligned}
f(x_1y_j)&=f(\sigma (x_1)y_{\tau(j)})=f(x_2y_{\tau(j)}).
\end{aligned}
\]
Similarly, $f(x_2y_j)=f(x_3y_{\tau(j)})$ and $f(x_3y_j)=f(x_1y_{\tau(j)})$. Consequently, for every $j$, 
\[
c_j=(a_j,b_j,d_j)=(f(x_2y_{\tau(j)}),f(x_3y_{\tau(j)}),f(x_1y_{\tau(j)}))=\rho(f(x_1y_{\tau(j)}),f(x_2y_{\tau(j)}),f(x_3y_{\tau(j)}))=\rho(c_{\tau(j)}).
\]

Since $\tau\in S_n$, we have
$S=\{c_1,...,c_n\}
=\{c_{\tau(1)},...,c_{\tau(n)}\}$.
So, 
\[
\rho(S)
=\{\rho(c_{\tau(1)}),...,\rho(c_{\tau(n)})\}
=\{c_1,...,c_n\}
=S,
\]
which contradicts our choice of $S$.
Similarly, $\sigma=(132)$ yields
$c_j=\rho^2(c_{\tau(j)})$
for every $j$, and hence
$S=\rho^2(S)$, which also contradicts our choice of $S$.
Therefore $\sigma=\mathrm{id}_{A_3}$.
\end{itemize}

\textbf{Step 4.} 
We prove $\tau=\mathrm{id}_{S_n}$. 

\begin{itemize}
    \item [] For every $j\in \{1,...,n\}$ and $i\in \{1,2,3\}$, preservation of $f$ gives $f(x_i y_j)=f(x_i y_{\tau(j)})$. 
Consequently, 
$c_j=(f(x_1 y_j),f(x_2 y_j), f(x_3 y_j))=(f(x_1 y_{\tau(j)}),f(x_2 y_{\tau(j)}),f(x_3 y_{\tau(j)}))=c_{\tau(j)}$.
Since all tuples 
$c_1,...,c_n$ in $S$ are pairwise distinct, we have $c_{j}=c_{\tau(j)}$ for all $j\in \{1,...,n\}$. Thus, $\tau=\mathrm{id}_{S_n}$.
\end{itemize}

Thus, $(\sigma,\tau)=(\mathrm{id}_{A_3},\mathrm{id}_{S_n})=\mathrm{id}_{A_3\times S_n}$ and so the only element of $A_3\times S_n$
preserving $f$ is the identity. So $f$ is a
distinguishing edge labeling using $q$ labels. Thus
$D'_{A_3\times S_n}\bigl(E(K_{3,n})\bigr)
\leq \left\lceil (n+1)^{1/3}\right\rceil$.
\end{proof}
Claims \ref{Claim 4.8} and \ref{Claim 4.9} together concludes the proof of Theorem \ref{Theorem 4.7}. 
\end{proof}


\subsection{Free actions}

Aliyu et al. \cite{AWTG2026} proved the following result.

\begin{lem}(\cite[Lemma 3.1]{AWTG2026})\label{Lemma 4.10}
{\em If $\Gamma$ acts freely on a set $V$, then $D_{\Gamma}(V) \leq 2$.}    
\end{lem}

\begin{prop}\label{Proposition 4.11}
{\em Let a group $\Gamma$ act freely on a graph $G=(V,E)$.
The following hold:
\begin{enumerate}
    \item $D''_{\Gamma}(V\cup E)\le 2$;
    \item If no nontrivial element of $\Gamma$ of order $2$
    swaps the endpoints of any edge, then $D'_{\Gamma}(E)\le 2$.
    \item The bounds in (1) and (2) are sharp.
\end{enumerate}
}
\end{prop}

\begin{proof}
(1). Since the action of $\Gamma$ on $V$ is free, $D_{\Gamma}(V)\leq 2$ by Lemma \ref{Lemma 4.10}. 
By Proposition \ref{Proposition 4.1}(2), we obtain
$D''_\Gamma(V\cup E)\le D_\Gamma(V)\le 2$.

(2). First, we prove the following claim.

\begin{claim}\label{Claim 4.12}
{\em For every edge $x=\{u,v\}\in E$, $\operatorname{Stab}_{\Gamma}(x)
=\{e\}$, where $e$ is the identity of $\Gamma$.}
\end{claim}

\begin{proof}
If $g\in\operatorname{Stab}_{\Gamma}(x)$ then
$
g\{u,v\}=\{g u,g v\}=\{u,v\}.
$
So, either $g u=u$ and $g v=v$, or
$g u=v$ and $g v=u$.
In the first case, since
the action is free, $g=e$. In the second
case,
$g^2 u=g(g u)=g v=u$, and so $g^2=e$.
Moreover, $g\ne e$. Hence $g$ has order $2$, contradicting the
assumption in (2). So,
$g=e$.
\end{proof}

By the orbit-stabilizer theorem, the orbit $\mathrm{Orb}_{\Gamma}(x)$ of each edge $x \in E$ has size
    $
    |\Gamma|/{|\mathrm{Stab}_{\Gamma}(x)|} = |\Gamma|,
    $
and the orbits partition $E$ into disjoint subsets of size $|\Gamma|$.
We define an edge labeling $\ell:E\rightarrow \{1,2\}$ of $G$. For each orbit $\mathrm{Orb}_{\Gamma}(x)$ of the edge $x \in E$, choose $z_x \in \mathrm{Orb}_{\Gamma}(x)$ and assign
    $\ell(z_x) = 1$. Assign $\ell(y) = 2$  for all other edges $y \in \mathrm{Orb}_{\Gamma}(x) \setminus \{ z_x \}$.
Since the orbits partition $E$, all edges are labeled. 
Let $g \in \Gamma$ preserve $\ell$, i.e.,
$
\ell(g y) = \ell(y)$ for all $y\in E$. 
Then $g$ maps edges labeled $1$ to edges labeled $1$.
Since every group element maps each orbit to itself, for every
$x\in E$ we have $gz_x\in\operatorname{Orb}_\Gamma(x)$.
Since each orbit $\mathrm{Orb}_{\Gamma}(x)$ contains a unique edge $z_{x}$ with the label $1$, it follows that $g z_x = z_x$ for every $x\in E$. By Claim \ref{Claim 4.12},
$
g\in\operatorname{Stab}_{\Gamma}(z_x)
=\{e\}.
$
Hence $g=e$.
So, $\ell$ is a distinguishing edge labeling of $G$ with respect to the action of $\Gamma$. Thus, $D'_{\Gamma}(E)\le 2$.

(3). Consider a cycle $C_n=(V(C_n), E(C_n))$, where $n \geq 3$, under the free action of $\mathbb{Z}_n$.  Since $|\mathbb{Z}_n| > 1$, Proposition \ref{Proposition 4.2} implies $D''_{\mathbb{Z}_n}(V(C_n) \cup E(C_n)) \geq 2$. Thus, by (1),
$D''_{\mathbb{Z}_n}(V(C_n) \cup E(C_n)) = 2$.
Consider the same action restricted to odd $n \geq 3$. Then $\mathbb{Z}_n$ has no non-identity element of order $2$. Similar reasoning yields $D'_{\mathbb{Z}_n}(E(C_n)) \geq 2$, which shows that the bound in $(2)$ is also sharp.
\end{proof}
\subsection{$n$-transitive actions}


\begin{lem}(\cite[Theorem 7.6A]{DM1996})\label{Lemma 4.13}
{\em Let $G$ be a finite group acting faithfully and exactly
$n$-transitively on a finite set $X$ of size $m\geq n$, where $n\geq4$.
Then, up to isomorphism, the tuple $(G,m,n)$ is one of
the combinations in Table 1, where $M_{11}$, $M_{12}$, $M_{23}$, and $M_{24}$ denote the
Mathieu groups:
}
\end{lem}

\begin{table}[h!]
\centering
\[
\begin{array}{c|c|c}
G&m&n\\ \hline
S_m&m&m\\
A_m&m&m-2\\
M_{11}&11&4\\
M_{12}&12&5\\
M_{23}&23&4\\
M_{24}&24&5.
\end{array}
\]
\caption{\em All possible combinations of $(G,m,n)$.}
\end{table}

\begin{thm}\label{Theorem 4.14}
{\em Let $G$ be a group acting faithfully and exactly $n$-transitively on a set $X$ of size $m\geq n$, where $n\geq4$.
Then
$D_G(X)\leq n+1$.}
\end{thm}

\begin{proof}
By Lemma~\ref{Lemma 4.13}, it suffices to
consider the six possibilities for $G$.

\medskip
\noindent
\textbf{Case 1: $G=S_m$.}
Then, 
$
D_{S_m}(X)=m=n\leq n+1.
$


\textbf{Case 2: $G\in \{A_m, M_{11}, M_{12}\}$.}
The action of $A_m$ on $X$ is sharply $(m-2)$-transitive, thus $A_m$ has a base of size $m-2$.
By Lemma~\ref{Lemma 4.3},
$D_{A_m}(X)\leq (m-2)+1=m-1=n+1$.
The action of $M_{11}$ on $11$ points is sharply
$4$-transitive and the action of $M_{12}$ on
$12$ points is sharply $5$-transitive. So, these
actions have bases of sizes $4$ and $5$, respectively.
By Lemma~\ref{Lemma 4.3}, we obtain
$
D_{M_{11}}(X)\leq 5=4+1=n+1
$
and
$
D_{M_{12}}(X)\leq 6=5+1=n+1.
$

\textbf{Case 3: $G=M_{23}$ or $G=M_{24}$.}
For $g\in G$, let $c(g)$ denote the number of cycles in the disjoint
cycle decomposition of $g$ in its permutation representation,
including cycles of length one.
Let $k=n+1$, and let $\mathcal{C}$ be the set of all $k$-colorings
of $X$. Since $|X|=m$, we have $|\mathcal{C}|=k^m$.
A $k$-coloring of $X$ is preserved by $g$ if and only if it is
constant on every cycle of $g$. Hence, the number of $k$-colorings
preserved by $g$ is $k^{c(g)}$.
Thus, the number of $k$-colorings preserved by at least one
nonidentity element of $G$ is
\[
\left\vert\bigcup_{1\ne g\in G}
\{f\in\mathcal{C}:g\text{ preserves }f\}\right\vert\leq \sum_{1\ne g\in G}k^{c(g)}.
\]
A $k$-coloring $f$ of $X$ is distinguishing if the identity element is the only element of $G$ that preserves $f$.
If
\begin{equation}\label{eq:union-bound}
\sum_{1\neq g\in G} k^{c(g)}<k^m
\end{equation}
then there exists at least one coloring preserved by no nonidentity
element of $G$, and hence a distinguishing $k$-coloring exists.
We recall the following datas from \cite{CCNPW85} (see Table \ref{Table 2}):
\begin{itemize}
    \item In $M_{23}$, there are $3795$ involutions, all of
    cycle type $1^7 2^8$. Every other nonidentity element has
    at most $11$ cycles.

    \item In $M_{24}$, there are $11385$ involutions of cycle
    type $1^8 2^8$ and $31878$ involutions of cycle type
    $2^{12}$. Every nonidentity element other than those of
    cycle type $1^8 2^8$ has at most $12$ cycles.
\end{itemize}

\begin{table}[h!]
\centering
\[
\begin{array}{c|r|r|c|c|c}
G & |G| & \text{Number of involutions}
& \multicolumn{2}{c|}{\text{Involution cycle types}}
& \max\limits_{|g|>2} c(g)\\
\cline{4-5}
& & & \text{Cycle type} & c(g) &\\
\hline
M_{23} & 10\,200\,960 & 3\,795
& 1^7 2^8 & 15 & 11\\
M_{24} & 244\,823\,040 & 11385
& 1^8 2^8 & 16 & 12\\
& & 
31878& 2^{12} & 12 & 
\end{array}
\]
\caption{\em Tabular representations of $M_{23}$ and $M_{24}$.}\label{Table 2}
\end{table}

First suppose that $G=M_{23}$. Since $n=4$, we take $k=5$. Then
\[
\begin{aligned}
\sum_{1\neq g\in M_{23}}5^{c(g)}
&\leq 3795\cdot 5^{15}
 +(10\,200\,960-3795)\cdot 5^{11}<5^{23}.
\end{aligned}
\]
Thus \eqref{eq:union-bound} holds, and consequently
$D_{M_{23}}(X)\leq 5=n+1$.
Now suppose that $G=M_{24}$. Here $n=5$, so we take $k=6$.
Since the elements of cycle type $1^8 2^8$ have $16$ cycles,
while every other nonidentity element has at most $12$ cycles,
we have, using $|M_{24}|=244\,823\,040$,

\[
\begin{aligned}
\sum_{1\ne g\in M_{24}}6^{c(g)}
&\le 11385\cdot 6^{16}
+31878\cdot 6^{12}+(244823040-(11385+31878))\cdot 6^{12}<6^{24}.
\end{aligned}
\]

Again, \eqref{eq:union-bound} holds, and hence
$
D_{M_{24}}(X)\leq 6=n+1.
$ 
This concludes the proof.
\end{proof}

\section{Actions on Cartesian powers of graphs}

Let $\Gamma$ be a group acting on a graph $G=(V,E)$.
The coordinatewise action of $\Gamma$ on the Cartesian power $V^k$ given by
$
g(v_1,...,v_k)=(gv_1,...,gv_k)
$
for all $g\in \Gamma$ and $(v_1,..., v_k)\in V^k$ induces an action of $\Gamma$ on the Cartesian power $G^k=(V^k,E_k)$ of $G$.

\vspace{-2mm}

\subsection{Hamming graphs}
For integers $k,n\ge 1$, the Hamming graph $H(k,n)=(V^k,E_k)$ is the Cartesian power $(K_n)^k$.
\begin{prop}\label{Proposition 5.1}
{\em Let $n\ge 3$, and let $S_n$ act naturally on $K_n$ and coordinatewise on the Hamming graph $H(k,n)=(K_n)^k=(V^k,E_k)$ for any integer $k\ge 2$. Then $D''_{S_n}(V^k\cup E_k)=D'_{S_n}(E_k)=2$.
}
\end{prop}

\begin{proof}
Let $[n]=\{1,...,n\}$ be the vertex set of $K_n$. Suppose $k\ge2$ and $n\geq 3$. 
For $1\le i<n$, let
$e_i=\{(i,1,...,1),(i+1,1,...,1)\}\in E_k$. Let $T=\{e_1,...,e_{n-1}\}$.
Define an edge labeling $\ell:E_k\to\{1,2\}$ of $H(k,n)$ using two labels by
\[
\ell(e)=
\begin{cases}
2, & \text{if } e\in T,\\
1, & \text{otherwise}.
\end{cases}
\]
Suppose $g\in S_n$ preserves $\ell$. Since
$T$ is the set of edges
labeled $2$,
$g(T)=T$.
In particular, for every $1\le i<n$,
$g(e_i)=\{(g(i),g(1),...,g(1)),
(g(i+1),g(1),...,g(1))\}\in T$. Since $k\ge2$,
the second coordinate of every endpoint of an edge in
$T$ is $1$. Hence $g(1)=1$.
Consequently,
$
g(e_i)=
\{(g(i),1,...,1),(g(i+1),1,...,1)\}
\in T$ for every
$1\le i<n$. Thus,
$
\{g(i),g(i+1)\}=\{j,j+1\}
$
for some $1\le j<n$. Therefore, $g$ induces an
automorphism of the path $P_n$ with vertex set $[n]$
and edge set $\{\{ i,i+1\}:1\le i<n\}$.
Since $g(1)=1$, it fixes the unique neighbor of $1$,
namely $2$. Inductively, it follows that $g(i)=i$
for every $i\in[n]$. Hence $g$ is the identity of
$S_n$, and $\ell$ is a distinguishing edge labeling.
Therefore,
$
D'_{S_n}(E_k)\le2.
$
Since $n\ge3$, we have $|S_n|>1$, and the
coordinatewise induced action of $S_n$ on $E_k$ is nontrivial. So $D'_{S_n}(E_k)\ge2$. Hence
$D'_{S_n}(E_k)=2$.
Since $H(k,n)$ is a connected graph of order at least $3$ for $k\geq 2$, $n\geq 3$, we obtain $D''_{S_n}(V^k\cup E_k)\le D'_{S_n}(E_k)=2$ by Proposition \ref{Proposition 4.2}(3).
Since $n\ge 3$, we have $|S_n|>1$, and the action of $S_n$ on $H(k,n)$ is nontrivial. By Proposition \ref{Proposition 4.2}, 
$D''_{S_n}(V^k\cup E_k)= 2$.
\end{proof}
\begin{remark}\label{Remark 5.2}
Mirafzal et al.~\cite{MZ2021} and Praeger et al.~\cite{PS2018} proved that \(\operatorname{Aut}(H(k,n))\cong S_n\wr S_k\). Imrich et al.~\cite[Theorem 1.1]{IK2006} proved that if \(G\neq K_2,K_3\) is a connected graph and \(k\geq 2\), then \(D(G^k)=2\). Thus,
    $D(H(k,n))=D_{\operatorname{Aut}(H(k,n))}(V^k)
    =D_{S_n\wr S_k}(V^k)=2$
for \(k\geq 2\) and \(n\geq 3\). Since $S_n\leq S_n\wr S_k$, it follows that \(D_{S_n}(V^k)\leq 2\). By Propositions~\ref{Proposition 4.1} and~\ref{Proposition 4.2},
$D''_{S_n}(V^k\cup E_k)=2$.
Pil\'{s}niak \cite[Theorem 3.2]{Pil2017} proved that every traceable graph of order at least \(7\) has distinguishing index at most \(2\). Since \(H(k,n)\) is traceable and has order at least $7$ for $k\geq 2$, $n\geq 3$,
$D'_{S_n}(E(H(k,n)))\leq D'_{Aut(H(k,n))}(E_k)=2$ by Proposition~\ref{Proposition 4.5}(1).
However, in Proposition~\ref{Proposition 5.1}, we provide an alternative self-contained proof of the result.
\end{remark}

\subsection{Coordinatewise action on Cartesian powers}

In \cite{Ali2026}, Aliyu proved the following.

\begin{thm}(\cite[Theorem 1]{Ali2026})
\label{Theorem 5.3}
{\em If $\Gamma$ acts transitively on $X$, then $D_{\Gamma}(X^k) \leq D_{\Gamma}(X) + k -1$ for all $k \geq 1$.}
\end{thm}

In the following theorem, we obtain an upper bound for the total distinguishing number of Cartesian powers of graphs that eliminates the need for the additional $k-1$ labels.

\begin{thm}\label{Theorem 5.4}
{\em Let a group $\Gamma$ act on a graph $G=(V,E)$. The following holds for all $k\geq 1$:
\begin{enumerate}
    \item If $D''_{\Gamma}(V^k\cup E_k)$ is the total distinguishing number of $G^k$ under the induced action of $\Gamma$ on $G^k$ and $D''_{\Gamma}(V\cup E)$ is the total distinguishing number of $G$ under the action of $\Gamma$ on $G$, then
$D''_{\Gamma}(V^k\cup E_k)\le D''_{\Gamma}(V\cup E)$.
    \item The bound in $(1)$ is sharp.
\end{enumerate}
}
\end{thm}

\begin{proof}
(1). Let $r=D''_{\Gamma}(V\cup E)$ and $\varphi:V\cup E\to\{1,...,r\}$ be a total distinguishing labeling of $G$.
We construct a total labeling
$
\Phi: V^k\cup E_k\to\{1,...,r\}
$
of $G^{k}$ as follows (cf. Fig. \ref{Figure 1}):

\begin{enumerate}
    \item[$(a)$] For a vertex $(v_1,...,v_k)\in V^k$, define
$\Phi(v_1,...,v_k)=\varphi(v_1)$.

    \item[$(b)$] In $G^k$, two vertices $v=(v_1,...,v_k)$ and $u=(u_1,...,u_k)$ are adjacent if they differ in exactly one coordinate, say $i$, and $\{v_i,u_i\}\in E$. 
We define

\[
\Phi\big(\{v,u\})=
\begin{cases}
\varphi(\{v_i,u_i\}), & \text{if } i=1,\\[4pt]
1, & \text{if } i>1.
\end{cases}
\]
\end{enumerate}

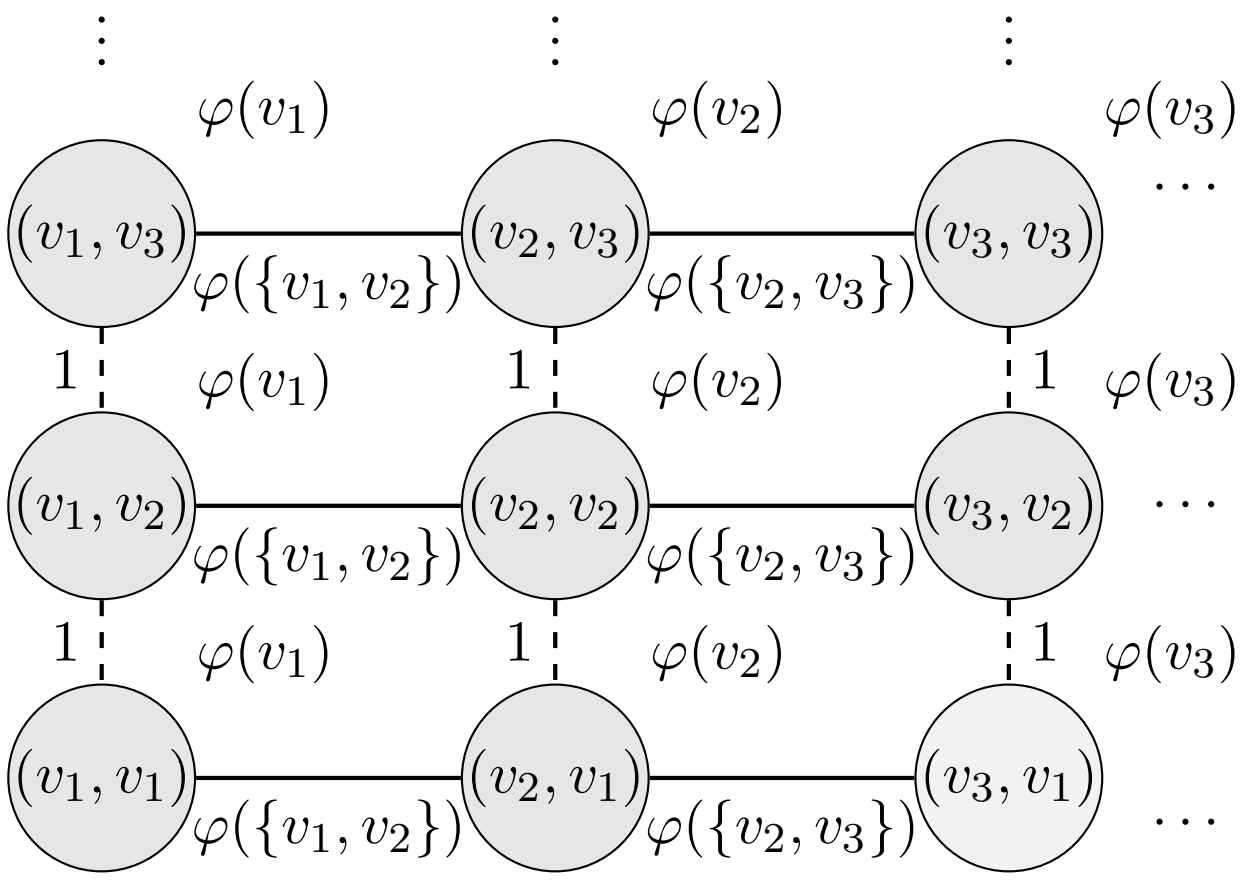
\begin{figure}[h]
\centering
\begin{tikzpicture}[scale=3,
    vertex/.style={circle, draw, inner sep=0pt, minimum size=16pt},
    edge/.style={thick}
]

\node[vertex, fill=gray!20] (a1) at (0,1.6) {$(v_1,v_3)$};
\node[above right=2pt of a1] {$\varphi(v_1)$};

\node[vertex, fill=gray!20] (a2) at (1,1.6) {$(v_2,v_3)$};
\node[above right=2pt of a2] {$\varphi(v_2)$};

\node[vertex, fill=gray!20] (a3) at (2,1.6) {$(v_3,v_3)$};
\node[above right=2pt of a3] {$\varphi(v_3)$};

\node[vertex, fill=gray!20] (b1) at (0,1) {$(v_1,v_2)$};
\node[above right=2pt of b1] {$\varphi(v_1)$};

\node[vertex, fill=gray!20] (b2) at (1,1) {$(v_2,v_2)$};
\node[above right=2pt of b2] {$\varphi(v_2)$};

\node[vertex, fill=gray!20] (b3) at (2,1) {$(v_3,v_2)$};
\node[above right=2pt of b3] {$\varphi(v_3)$};

\node[vertex, fill=gray!20] (c1) at (0,0.4) {$(v_1,v_1)$};
\node[above right=2pt of c1] {$\varphi(v_1)$};

\node[vertex, fill=gray!20] (c2) at (1,0.4) {$(v_2,v_1)$};
\node[above right=2pt of c2] {$\varphi(v_2)$};

\node[vertex, fill=gray!10] (c3) at (2,0.4) {$(v_3,v_1)$};
\node[above right=2pt of c3] {$\varphi(v_3)$};

\draw[edge] (a1) -- node[below] {$\varphi(\{v_1,v_2\})$} (a2);
\draw[edge] (a2) -- node[below] {$\varphi(\{v_2,v_3\})$} (a3);

\draw[edge] (b1) -- node[below] {$\varphi(\{v_1,v_2\})$} (b2);
\draw[edge] (b2) -- node[below] {$\varphi(\{v_2,v_3\})$} (b3);

\draw[edge] (c1) -- node[below] {$\varphi(\{v_1,v_2\})$} (c2);
\draw[edge] (c2) -- node[below] {$\varphi(\{v_2,v_3\})$} (c3);

\draw[edge, dashed] (a1) -- node[left] {$1$} (b1);
\draw[edge, dashed] (b1) -- node[left] {$1$} (c1);

\draw[edge, dashed] (a2) -- node[left] {$1$} (b2);
\draw[edge, dashed] (b2) -- node[left] {$1$} (c2);

\draw[edge, dashed] (a3) -- node[right] {$1$} (b3);
\draw[edge, dashed] (b3) -- node[right] {$1$} (c3);

\node at (0,2.06) {$\vdots$};
\node at (1,2.06) {$\vdots$};
\node at (2,2.06) {$\vdots$};

\node at (2.4,0.3) {$\cdots$};
\node at (2.4,1) {$\cdots$};
\node at (2.4,1.7) {$\cdots$};

\end{tikzpicture}

\caption{\em Total labeling of the Cartesian power $P_n^2$ of a path $P_n=v_1v_2v_3...$.}
\label{Figure 1}
\end{figure}

It suffices to show that $\Phi$ is a distinguishing total labeling of $G^k=(V^k,E_k)$ with respect to the action of $\Gamma$.
Suppose $g\in \Gamma$ preserves $\Phi$.
We have the following:

\begin{enumerate}
    \item[$(i)$] Fix $v_1\in V$. For any $(v_1,...,v_k)\in V^k$, since $g$ preserves $\Phi$,
\[
\varphi(gv_1)
=\Phi(gv_1,...,gv_k)
=\Phi(g(v_1,...,v_k))
=\Phi(v_1,...,v_k)
=\varphi(v_1).
\]
    \item[$(ii)$] Let $\{v_1, u_1\} \in E$. Consider the vertices 
$
v = (v_1, v_2, ..., v_k)$ and $u = (u_1, v_2, ..., v_k)
$
in $V^k$ for some $v_2, ..., v_k \in V$. Then $\{v, u\} \in E_k$. By the construction of $\Phi$, we obtain \[\Phi(\{v, u\}) = \varphi(\{v_1, u_1\}).\]
Moreover, $\{g v,g u\}\in E_k$, and
$
\Phi(\{gv,g u\})=\varphi(\{gv_1,gu_1\}).
$
As $g$ preserves $\Phi$, 
\[
\varphi(g\{v_1,u_1\})=\varphi(\{gv_1,gu_1\})=\Phi(\{gv,gu\})=\Phi(g\{v,u\})=\Phi(\{v,u\})=\varphi(\{v_1,u_1\}).
\]
\end{enumerate}

By $(i)$ and $(ii)$, $g$ preserves the total labeling $\varphi$. Since $\varphi$ is a total distinguishing labeling, it follows that
$g\in \mathrm{Stab}_{\Gamma}(V\cup E)=\mathrm{Stab}_{\Gamma}(V)$ by Proposition \ref{Proposition 4.1}(1).
Since $\Gamma$ acts coordinatewise on $V^k$, we also have
$\operatorname{Stab}_{\Gamma}(V^k)
=\operatorname{Stab}_{\Gamma}(V)$.
Applying Proposition \ref{Proposition 4.1}(1) to $G^k$, we obtain
\[
g\in\mathrm{Stab}_{\Gamma}(V)=\mathrm{Stab}_{\Gamma}(V^k)=\mathrm{Stab}_{\Gamma}(V^k\cup E^k).
\]
Thus, $\Phi$ is a distinguishing total labeling of $G^k$.

(2). Consider the faithful action of $\mathbb{Z}_2$ on the path $P_n=(V, E)$ for $n\geq 3$. 
Since $|\mathbb{Z}_2|>1$, the action is nontrivial. 
Proposition \ref{Proposition 4.2} implies $D''_{\mathbb{Z}_2}(V \cup E) \geq 2$.  A  total labeling of $P_n$ with $2$ labels
assigning $2$ to an end-vertex of $P_n$ and $1$ to all other vertices and edges is distinguishing, which shows $D''_{\mathbb{Z}_2}(V\cup E)=2$. 
For any $k\geq 1$, since the coordinatewise action of $\mathbb{Z}_2$ on $P_{n}^{k}=(V^k,E_k)$ is nontrivial, Proposition \ref{Proposition 4.2} implies $D''_{\mathbb{Z}_2}(V^k \cup E_k) \geq 2$. By (1), $D''_{\mathbb{Z}_2}(V^k \cup E_k) = 2$.
\end{proof}
\subsection{Wreath product action on Cartesian powers}

Aliyu \cite{Ali2026} obtained the following result.

\begin{thm}(\cite[Theorem 3]{Ali2026})
\label{Theorem 5.5}
{\em Let $\Gamma$ act faithfully on a finite set $X$, and let $W = \Gamma \wr S_m$ act in the natural product action on $X^m$. Then
$D_{W}(X^m) \leq D_{\Gamma}(X)+m-1.$}
\end{thm}

We prove that a similar upper bound holds for the total distinguishing number of the Cartesian power graph under the wreath product action.

\begin{thm}\label{Theorem 5.6}
{\em Let $m\geq 2$ and $\Gamma$ be a group acting faithfully on a graph $G=(V,E)$ where $|V|\geq 2, |E|\geq 1$, and $|\Gamma|>1$.
Then the action of $W=\Gamma\wr S_m$ on $V^m$ given by
\[
(g_1,...,g_m;\sigma)(v_1,...,v_m)
=
(g_1v_{\sigma^{-1}(1)},...,g_mv_{\sigma^{-1}(m)})
\]
for all $g_j \in \Gamma$ for $1 \le j \le m$ and $(v_1, ..., v_m) \in V^m$, induces an action of $W$ on the Cartesian power $G^{m}=(V^{m},E_m)$.
If $D''_{W}(V^m\cup E_m)$ is the total distinguishing number of $G^m$ under the induced action of $W$ on $G^m$ and $D''_{\Gamma}(V\cup E)$ is the total distinguishing number of $G$ under the action of $\Gamma$ on $G$,
then
$
D''_{W}(V^m\cup E_m)\le D''_{\Gamma}(V\cup E)+m-1.
$}
\end{thm}

\begin{proof}
Let $r=D''_{\Gamma}(V\cup E)$.
Since $|\Gamma|>1$ and the action of $\Gamma$ on $G$ is
faithful, we have $r\ge2$ by Proposition \ref{Proposition 4.2}.
Since $E\neq\emptyset$ and
$r\ge2$, we may choose a total distinguishing labeling
$\phi:V\cup E\to\{1,...,r\}$ such that for some edge
$\{u,w\}\in E$, $\phi(\{u,w\})\ne1$.\footnote{If every edge has label $1$ with respect to $\phi$, relabel one edge with
label $2$. The new labeling $\psi$ remains total distinguishing;
any element of $\Gamma$ preserving $\psi$ must fix the relabeled edge, since it is the unique edge with label
$2$. It also preserves
$\phi$, and hence belong to the pointwise stabilizer of $V\cup E$ which is the identity of $\Gamma$ since the action of $\Gamma$ on $G$ is faithful.}
Define a labeling
$\Phi:V^m\cup E_m\to\{1,...,r+m-1\}$ as follows (cf. Fig. \ref{Figure 2}):

\begin{enumerate}
    \item For $v=(v_1,...,v_m)\in V^m$, 

\begin{enumerate}
\item If $v_1=v_2=\cdots=v_m$, then define $\Phi(v)=\phi(v_1)$.

\item Otherwise, if
$i=\min\{j\ge2 : v_j\ne v_1\}$, then let
$
\Phi(v)=r+i-1.
$
\end{enumerate}
   \item Let
$v=(v_1,...,v_m)$ and
$u=(v_1,...,v_m)$
be adjacent vertices in $G^m$ with $\{v_i,u_i\}\in E$ for some $1\leq i\leq m$.
Define
\[
\Phi(\{v,u\})=
\begin{cases}
\phi(\{v_i,u_i\}), & i=1,\\
1, & i\ge2.
\end{cases}
\]
\end{enumerate}

\begin{figure}[h]
\centering
\begin{tikzpicture}[scale=2.5,
    vertex/.style={circle, draw, inner sep=0pt, minimum size=16pt},
    diag/.style={circle, draw, fill=gray!20, inner sep=0pt, minimum size=16pt},
    edge/.style={thick}
]

\node[diag] (a1) at (0,1.7) {$(v_1,v_1)$};
\node[vertex] (a2) at (1.5,1.7) {$(v_2,v_1)$};
\node[vertex] (a3) at (3,1.7) {$(v_3,v_1)$};

\node[vertex] (b1) at (0,1) {$(v_1,v_2)$};
\node[diag] (b2) at (1.5,1) {$(v_2,v_2)$};
\node[vertex] (b3) at (3,1) {$(v_3,v_2)$};

\node[vertex] (c1) at (0,0.3) {$(v_1,v_3)$};
\node[vertex] (c2) at (1.5,0.3) {$(v_2,v_3)$};
\node[diag] (c3) at (3,0.3) {$(v_3,v_3)$};

\foreach \row in {a,b,c} {
    \draw[edge] (\row1) -- (\row2);
    \draw[edge] (\row2) -- (\row3);
}

\foreach \col in {1,2,3} {
    \draw[edge, dashed] (a\col) -- (b\col);
    \draw[edge, dashed] (b\col) -- (c\col);
}

\node[above right=2pt of a1] {$\phi(v_1)$};
\node[above right=2pt of b2] {$\phi(v_2)$};
\node[above right=2pt of c3] {$\phi(v_3)$};

\node[above right=2pt of a2] {$r+1$};
\node[above right=2pt of a3] {$r+1$};
\node[above right=2pt of b1] {$r+1$};
\node[above right=2pt of b3] {$r+1$};
\node[above right=2pt of c1] {$r+1$};
\node[above right=2pt of c2] {$r+1$};

\draw[edge] (a1) -- node[below] {$\phi(\{v_1,v_2\})$} (a2);
\draw[edge] (a2) -- node[below] {$\phi(\{v_2,v_3\})$} (a3);

\draw[edge] (b1) -- node[below] {$\phi(\{v_1,v_2\})$} (b2);
\draw[edge] (b2) -- node[below] {$\phi(\{v_2,v_3\})$} (b3);

\draw[edge] (c1) -- node[below] {$\phi(\{v_1,v_2\})$} (c2);
\draw[edge] (c2) -- node[below] {$\phi(\{v_2,v_3\})$} (c3);

\draw[edge, dashed] (a1) -- node[left] {$1$} (b1);
\draw[edge, dashed] (b1) -- node[left] {$1$} (c1);

\draw[edge, dashed] (a2) -- node[left] {$1$} (b2);
\draw[edge, dashed] (b2) -- node[left] {$1$} (c2);

\draw[edge, dashed] (a3) -- node[right] {$1$} (b3);
\draw[edge, dashed] (b3) -- node[right] {$1$} (c3);

\node at (0,2.06) {$\vdots$};
\node at (1.5,2.06) {$\vdots$};
\node at (3,2.06) {$\vdots$};

\node at (3.4,0.3) {$\cdots$};
\node at (3.4,1) {$\cdots$};
\node at (3.4,1.7) {$\cdots$};

\end{tikzpicture}
\caption{\em Total labeling of $P_n^2$ under the natural action of $\Gamma \wr S_2$.}
\label{Figure 2}
\end{figure}
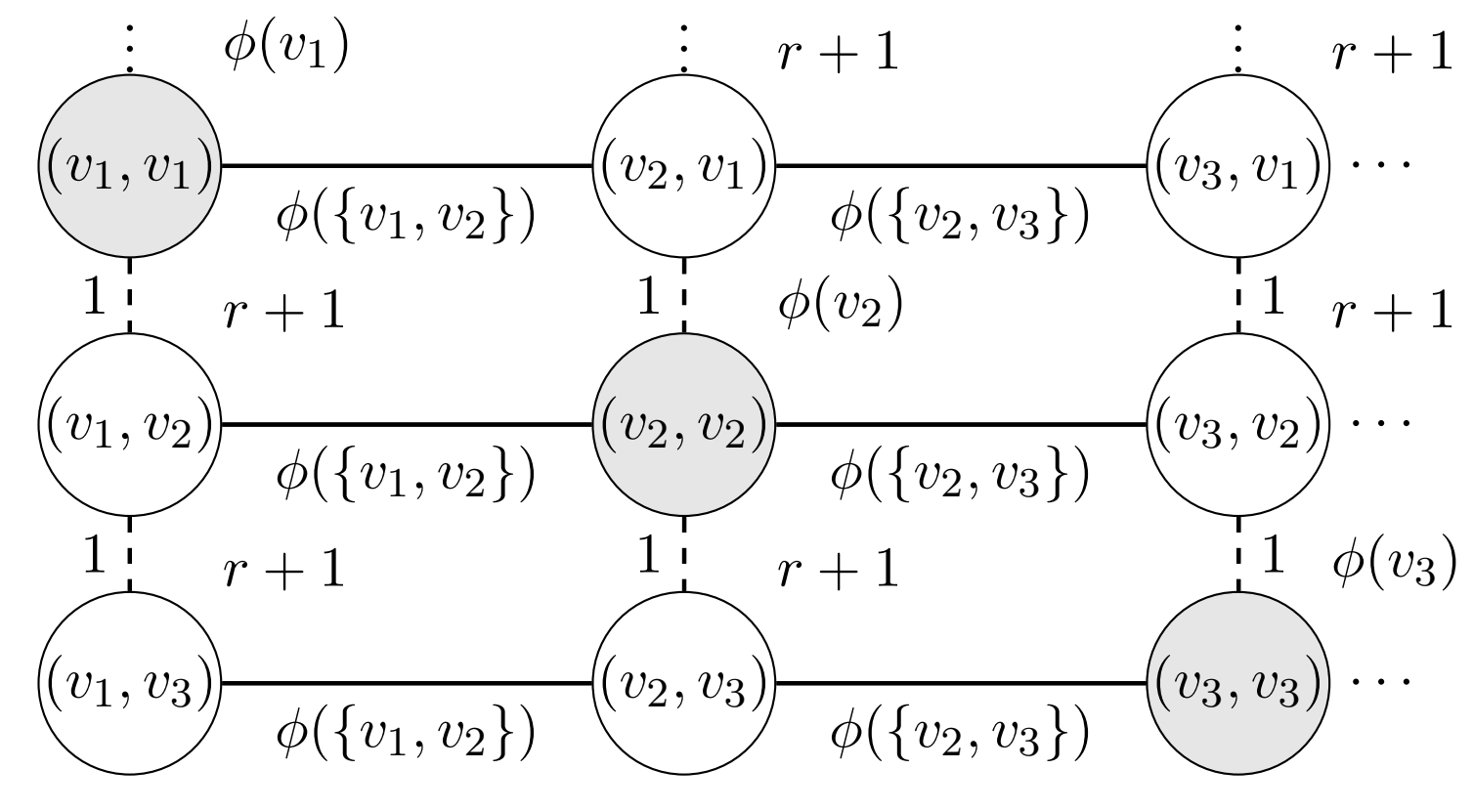

The following lemma is immediate from the definition of $\Phi$.

\begin{lem}\label{Lemma 5.7}
{\em For $v=(v_1,...,v_m)\in V^m$, $\Phi(v)\in\{1,...,r\}$ if and only if $v_1=...=v_m$.}
\end{lem}

We show that $\Phi$ is distinguishing.
Let $(g_1,...,g_m;\sigma)\in W$ preserve $\Phi$. 
Then
\[
\Phi((g_1,...,g_m;\sigma) u)=\Phi(u)
\quad \text{for all } u\in V^m\cup E_m.
\]

\textbf{Step 1.} We prove  $g_1=g_2=\cdots=g_m$.

\begin{itemize}
    \item [] 
For any $a\in V$, consider $v=(a,...,a)$.
Then $\Phi((g_1,...,g_m;\sigma) v)=\Phi(v)$. So, $\Phi(g_1a,...,g_ma)=\phi(a)\in\{1,...,r\}$.  
By Lemma \ref{Lemma 5.7}, we get
$
g_1a=g_2a=\cdots=g_ma.
$
Thus, $g_i a=g_j a$ for all $a\in V$ and all $i,j$. Consequently,
$g_j^{-1}g_i$ fixes every element of $V$. 
Since the action of $\Gamma$ 
on $V$ is faithful, $g_j^{-1}g_i=e$, and so
$
g_1=\cdots=g_m=g
$
for some $g\in \Gamma$.
\end{itemize}

\medskip

\textbf{Step 2.} We prove that $\sigma=\mathrm{id}_{S_m}$, where $\mathrm{id}_{S_m}$ is the identity of $S_m$.

\begin{itemize}
    \item []
For the sake of contradiction, suppose $\sigma\ne\mathrm{id}_{S_m}$. Either there exists  a
minimum 
$j\in\{2,...,m\}$ with $\sigma(j)\ne j$ and
$\sigma(j)\ge2$, or $\sigma=(1\, j)$ is a transposition of $1$ and $j$ for some $j\geq 2$.\footnote{In the second case, every moved $j\ge2$ satisfies $\sigma(j)=1$.
Since $\sigma$ is a permutation, there is at most one such
$j$. Since $\sigma\ne\mathrm{id}$, there must be exactly one,
say $j$. Furthermore, $\sigma(1)=j$, and so $\sigma=(1\, j)$.}
In either cases, we pick $a,b\in V$ with $a\ne b$, and let
$v=(v_1,...,v_m)$, where $v_j=b$, and
$v_\ell=a$ for all $\ell\neq j$. 
Since 
$j=\min\{t\ge2 : v_t\ne v_1\}$, we have
$\Phi(v)=r+j-1$.
\end{itemize}

\textbf{Case 1: there exists a minimum
$j\in\{2,...,m\}$ such that $\sigma(j)\ne j$ and
$\sigma(j)\ge2$.}

\begin{itemize}
    \item []

Since the only coordinate of $v$ equal to $b$ is the $j$-th coordinate, 
the tuple $(v_{\sigma^{-1}(1)},...,v_{\sigma^{-1}(m)})$ has $b$ 
in position $\sigma(j)$. Let $i=\sigma(j)$.  
Since $a\ne b$ and $g$ acts bijectively on $V$, $ga\ne gb$. 
The first coordinate where $(gv_{\sigma^{-1}(1)},...,gv_{\sigma^{-1}(m)})$ differs from its first entry is $i$. So, 
\[ \Phi((g,...,g;\sigma)v)=\Phi(gv_{\sigma^{-1}(1)},...,gv_{\sigma^{-1}(m)})=r+i-1.\]
                                  
Since $i=\sigma(j)\ne j$, we obtain
$\Phi((g,...,g;\sigma) v)\ne \Phi(v)$,
contradicting the assumption that $(g,...,g;\sigma)$ preserves $\Phi$.
\end{itemize}

\textbf{Case 2: $\sigma=(1\,j)$ with $j\ge3$.}

\begin{itemize}
    \item []
The tuple $(gv_{\sigma^{-1}(1)},...,gv_{\sigma^{-1}(m)})$ has $gb$ in its first coordinate and $ga$
in all other coordinates. Since $ga\ne gb$, its first
coordinate differing from its first entry is at position two. Thus
$
\Phi((g,...,g;\sigma) v)=r+1\ne r+j-1,
$
again a contradiction.
\end{itemize}
\textbf{Case 3: $\sigma=(1\,2)$.}

\begin{itemize}
    \item []
We picked the
distinguishing labeling $\phi$ so that for some edge
$e=\{u,w\}\in E$, $\phi(e)\ne1$.
Consider the edge of $G^m$ joining
$x=(u,u,...,u)$ and $y=(w,u,...,u)$. 
Then $\Phi(\{x,y\})=\phi(e)\ne1$.
Under $(g,...,g;(1\,2))$, the edge $\{x,y\}$ maps to the edge
\[\{(gu,gu,...,gu), (gu,gw,gu,...,gu)\},\]
which differs only in coordinate $2$ and therefore has
$\Phi$-label $1$. This contradicts the assumption that
$(g,...,g;\sigma)$ preserves $\Phi$.

\end{itemize}

\medskip
\textbf{Step 3.} We prove that $g=e$, where $e$ is the identity of $\Gamma$.

\begin{itemize}
    \item []     
\begin{enumerate}[$(a)$.]
    \item Since $(g_1,...,g_m;\sigma)=(g,...,g;\mathrm{id}_{S_m})$ preserves $\Phi$, applying it to the constant tuple
    $(a,...,a)$ gives $\phi(ga)=\phi(a)$ for all  $a\in V$.
    \item For any edge $\{a,b\}\in E$, the edge joining $u=(a,...,a)$ and $v=(b,a,...,a)$ has label $\Phi(\{u,v\})=\phi(\{a,b\})$ by the definition of $\Phi$ on the edge set $E_m$.

    \item 
Now,
\begin{align*}
(g,...,g;\mathrm{id}_{S_m}) \{u,v\} &= \{(g,...,g;\mathrm{id}_{S_m}) u,(g,...,g;\mathrm{id}_{S_m}) v\} \\
&= \{(ga,...,ga),(gb,ga,...,ga)\}.
\end{align*}

Since $\Phi((g,...,g;\mathrm{id}_{S_m}) y)=\Phi(y)$ for all $y\in E_m$, 
$
\Phi(\{(ga,...,ga),(gb,ga,...,ga)\})$
=
$\Phi(\{u,v\})$, which is $\phi(a,b)$ by step $(b)$.
By the definition of $\Phi$ on the edge set $E_m$, 
\[
\Phi(\{(ga,...,ga),(gb,ga,...,ga)\})
=
\phi(\{ga,gb\}).
\]
Therefore, $\phi(\{ga,gb\})=\phi(\{a,b\})$
for all $\{a,b\}\in E$.

\end{enumerate}

\smallskip
By steps $(a)$ and $(c)$, $g$ preserves the total labeling $\phi$ on $V\cup E$.
Since $\phi$ is a total distinguishing labeling of $G$ and the action of $\Gamma$ on $G$ is faithful, it follows that $g=e$.
\end{itemize}

By Steps 1, 2 and 3, the only element of $W$ that preserves
$\Phi$ is $((e,...,e);\operatorname{id}_{S_m})$, the identity
element of $W$. Hence, $\Phi$ is a total distinguishing labeling
using at most $r+m-1$ labels.
\end{proof}

\end{document}